\documentclass[graphicx,amscd,amssymb,amsfonts,enumerate,verbatim, amsmath,amsthm,color,12pt,righttag]{amsart}
\usepackage{color}
\usepackage{txfonts}

 \newtheorem{theorem}{Theorem}[section]
 \newtheorem{Def}[theorem]{Definition}
 \newtheorem{Pro}[theorem]{Proposition}
 \newtheorem{Lem}[theorem]{Lemma}

 \newtheorem{Example}[theorem]{Example}

 \numberwithin{equation}{section}
 
\newcommand{\A}{{\mathcal A}}
\newcommand{\B}{{\mathcal B}}
\newcommand{\C}{{\mathcal C}}
\newcommand{\D}{{\mathcal D}}
\newcommand{\E}{{\mathcal E}}

\newcommand{\R}{{\mathbb R}}
\newcommand{\Z}{{\mathbb Z}}

\begin{document}

\title {Spectrality of Self-Similar Tiles}

%\date{}
\author{Xiaoye Fu}
\address{Department of Mathematics, The Chinese University of Hong Kong,  Hong Kong}
\author {Xinggang He}
\address {College of Mathematics and Statistics, Central China Normal University, Wuhan, China}
\author{Ka-Sing Lau}
\address{Department of Mathematics, The Chinese University of Hong Kong, Hong Kong}

\email{Xiaoyefu@gmail.com}
\email{xingganghe@163.com}
\email{kslau@math.cuhk.edu.hk}

%\thanks {}

%\date{\today}

\thanks{The research is supported in part by the HKRGC Grant, and the CNSF (no. 11171100).  The second author is also supported by the National Natural
Science Foundation of China 11271148. }

\maketitle

\begin{abstract}
We call a set $K \subset {\mathbb R}^s$ with positive Lebesgue
measure a {\it spectral set} if $L^2(K)$ admits an exponential
orthonormal basis. It was conjectured that $K$ is a spectral set if and only if $K$ is a tile (Fuglede's conjecture). Despite the conjecture was proved to be false on ${\mathbb R}^s$, $s\geq 3$ ([T], [KM2]), it still poses challenging questions with additional assumptions. In this paper, our additional assumption is self-similarity.  We study the spectral properties for the class of self-similar tiles $K$ in ${\mathbb R}$ that has a product structure on the associated digit sets. We show that any strict product-form tiles and the associated modulo product-form tiles are spectral sets.  As for the converse question, we give a pilot study for the self-similar set $K$ generated by arbitrary digit sets with four elements. We investigate the zeros of its Fourier transform due to the orthogonality, and verify Fuglede's conjecture for this special case.
\end{abstract}

%\tableofcontents

\medskip

\section{\bf Introduction}

\medskip

Let  $\mu$ be  a finite regular Borel measure on ${ \mathbb R}^s$ with compact support $K$.
 We say that $\mu$ is a {\it spectral measure} if $L^2(K, \mu)$ admits an exponential orthonormal basis, i.e., there exists a countable set
 $\Lambda \subset {\mathbb R}^s$ such that $\{e^{2\pi i \langle\lambda,  x\rangle} : \lambda \in \Lambda\}$ is an orthonormal
 basis for $L^2(K, \mu)$ (we always assume $0 \in \Lambda$). If, in addition, $K$ has positive Lebesgue measure and $\mu$ is the
 Lebesgue measure on $K$, then we say that $K$ is a {\it spectral set}. If ${\mathcal A} \subset {\mathbb R}^s$
is a finite set and $\mu= \delta_{\mathcal A}$ is the point mass measure, then it is easy to see that $\mu$ is a spectral measure
if and only if there exists a finite set $\Lambda \subset {\mathbb R}^s$ such that $[e^{2\pi i \langle \lambda, a\rangle}]_{a \in \A, \lambda \in
\Lambda}$
is an orthogonal matrix, we also say  $\A$ is a spectral set.

\medskip

The notion of spectral set was first introduced by Fuglede [F] in his study of the extendability of the commuting partial differential operators $\{\partial_j\}_{j=1}^s$  on  $C_c(K)$ to commuting self-adjoint operators on $L^2(K)$. He characterized that  $K$ has to be a spectral set. Furthermore, he showed that a lattice tile is a spectral set, and raised the far reaching conjecture that {\it a bounded Borel set $K \subset {\mathbb R}^s$ is a spectral set if and only $K$ is a tile}.

\medskip

The Fuglede's conjecture generated a lot of interest in the last four decades. There were many partial results, and involved a lot of techniques from harmonic analysis,  tiling theory, convex geometry and geometry of numbers ([P], [LW4], [IKP], [IKT1,2], [$\L$1,2]). The conjecture was finally disproved by Tao [T] that a spectral set is not necessary a tile for ${\mathbb R}^s, s\geq 5$, and subsequently improved by Kolountzakis and Matolcsi [KM1,2] that the conjecture is false in both directions for  $s\ge 3$. Nevertheless, the conjecture is still valid for ${\mathbb R}^s$, $s\leq 2$, with additional conditions (e.g., convexity). The problem also led to the study of the spectral property of the discrete integer tiles ([LW4], [PW]) and the Cantor-type singularly continuous measures ([JP], [$\L$W], [HuL], [D], [DHS], [DHL]).

\bigskip

In this paper, our main purpose is to initiate a study on the
spectral property of self-similar tiles on ${\mathbb R}$. Let $A$ be an $s\times s$
expanding matrix (i.e., all eigenvalues have moduli $>1$) with
integer entries and $|\det A|=b$. Let ${\mathcal D} \subset {\mathbb
Z}^s$ be a finite set and call it a {\it digit set}. It follows that
for the pair $(A, {\mathcal D})$, there exists a unique compact set
$K:= K(A, {\mathcal D})$ satisfying the set-valued equation
\begin {equation} \label{eq1.1}
AK = K+ {\mathcal D}.
\end{equation}
Alternatively, $K$ can be expressed as a set of radix expansion with base $A$, i.e., $K = \{\sum_{k=1}^\infty A^{-k}d_k: \ d_k \in {\mathcal D}\}$. The self-affine set $K(A, \D)$ also associates with a probability measure $\mu := \mu_{A, \D}$ ({\it self-affine measure}) [H] defined by
\begin {equation} \label {eq1.2}
\mu (\cdot)= \frac{1}{\# {\mathcal D}} \sum_{d\in {\mathcal D}}\mu(A \cdot -d).
\end{equation}
 It is well-known that if $|\det A| = \#{\mathcal D}=b$ and $K$ has non-void interior, then $K$ is a translational tile in ${\mathbb R}^s$ [B]. We call such $K$ a {\it self-affine tile} (or a {\it self-similar} tile if $A$ is a scalar multiple of an orthogonal matrix), and  call ${\mathcal D}$ a {\it tile digit set} with respect to $A$.  These tiles are referred to as {\it fractal tiles} because their boundaries are usually fractals. In such case, the corresponding self-affine measure $\mu$ in (\ref{eq1.2}) is the Lebesgue measure restricted on $K$. There are many literatures concerning such tiles (see [LW1,3], [LR],  [LLR1,2] and the references therein). One of the fundamental questions in the theory is:

{\it For a given $A$, determine ${\mathcal D}$ to be a tile digit set, i.e., $K(A, {\mathcal D})$ is a self-affine tile.}

\vspace {0.1cm}

\noindent Necessarily, a tile digit set ${\mathcal D}$ must satisfy $\#{\mathcal D} = |\det A| = b$. It is well-known that if ${\mathcal D}$ is a complete residue set modulo $A$ (i.e.,  $\D \equiv {\mathbb Z}^s/A{\mathbb Z}^s)$, then ${\mathcal D}$ is a tile digit set [B]; it is called a {\it standard tile digit set} in [LW1]. In fact, this characterizes the tile digit sets when $|\det A|= b$ is a prime (see [K] for ${\mathbb R}$, and [LW1] and [HL] for ${\mathbb R}^s$ with some additional conditions on ${\mathcal D}$).  When $|\det A| =b$ is not a prime, Lagarias and Wang [LW1] extended a direct sum setup of Odlyzko [O] to the class of {\it product-form} tile digit set:
\begin{equation} \label {eq1.3}
  {\mathcal D} = \E_0 \oplus A^{\ell_1} \E_1 \oplus \cdots \oplus A^{\ell_k} {\E_k},
\end{equation}
where $\E = \E_0 \oplus  \E_1 \cdots \oplus  {\E_k}$ is a complete residue set modulo $A$, and $0 \leq \ell_1\leq \dots \leq \ell_k$. This idea has been further investigated by Lai, Lau and Rao ([LR], [LLR1,2]). They extended the form to certain  {\it modulo product-form} and  {\it higher order modulo product-form} in $\R$, by taking appropriate modulo operation on each direct summand (see Section 3).  These classes include all the known tile digit sets, and it is speculated the higher order modulo product-forms may cover all tile digit sets in ${\mathbb R}$ [LLR1,2].

\medskip

Our main purpose is to study the spectral property of the modulo product-forms on ${\mathbb\mathbb\mathbb R}$. We assume that ${\mathcal D}\subset \Z^+, \ 0 \in {\mathcal D}$ and gcd$({\mathcal D})=1$. It is known that for the product-form in (\ref {eq1.3}), there is a finite set ${\mathcal A} \subset {\mathbb Z}$ such that $K(b, {\mathcal D}) =K(b, \E) \oplus \A$ [LW2], we prove (see Proposition \ref{th3.2}, Theorems \ref{th3.3}, \ref{th3.4}, \ref{th3.6}.)

\medskip

\begin {theorem} \label {th1.1}
 Let ${\mathcal D} \subset \Z^+$ be a product-form tile digit set with $\#{\mathcal D} = b$. Then $K(b, {\mathcal D})$ is a spectral set if and only if $\A$ is a spectral set. In particular, if $\D$ is a product-form with respect to \ $\E_0 \oplus  \E_1 \cdots \oplus  {\E_k}= {\Bbb  Z}_b$, then $K(b, \D)$ is spectral set.
\end{theorem}

\medskip

\begin{theorem}\label{th1.2}
Let ${\mathcal D}$ be a modulo product-form (as in  Definition \ref{th3.5} ) with respect to  $\E_0 \oplus  \E_1 \cdots \oplus  {\E_k} = {\mathbb Z}_b$. Then $K(b,{\mathcal D})$ is a spectral set, and any spectrum of the corresponding product-form $K(b,{\mathcal D}^{\prime})$  is also a spectrum of $K(b,{\mathcal D})$.
\end{theorem}

\medskip

 The first theorem allows us to reduce the spectrality of the  product-forms to the simpler discrete spectral sets. Note that for $b= p^\alpha$ or $pq$ where $p, q$ are distinct primes and $\alpha >0$ is an integer, the tile digit sets are characterized by the modulo product-forms defined by $\E_0 \oplus  \E_1 \cdots \oplus  {\E_k} = {\mathbb Z}_b$ ([LR], [LLR2]). Hence, it follows from Theorem \ref{th1.2} that

\medskip

\begin {Pro} \label {th1.3}
Let  $b = p^\alpha$ or $b=pq$, where $p, q$ are distinct primes, and let ${\mathcal D}$ be a tile digit set associated with $\#{\mathcal D} = b$. Then the self-similar tile $K(b, {\mathcal D})$ is a spectral set.
\end{Pro}

\medskip

In view of the Fuglede problem, we also consider the reverse question:  {\it suppose $\#{\mathcal D} = b$, if the self-similar measure $\mu_{b, {\mathcal D}}$ in (\ref{eq1.2}) is a spectral measure, will $K(b, {\mathcal D})$ be a self-affine tile?} This seems to be a more difficult question. In order to gain more insight,  we prove in detail the following very special case that ${\mathcal D}$ contains four elements.

\bigskip

\begin {theorem} \label {th1.4}
Suppose ${\mathcal D} \subset \Z^+$ and $\#{\mathcal D} =4$. Then $K(4, {\mathcal D})$ is a self-affine tile if and only if it is a spectral set.  In this case  ${\mathcal D} = \{0, a, 2^t\ell, a + 2^t \ell'\}$ where $a, t, \ell, \ell'$ are odd integers, and $K(4, {\mathcal D})$  has a spectrum  given by
$\Lambda= \Z +\sum_{j=1}^k\frac{1}{4^j}\{0,1\}$ with $k = (t-1)/2$.
\end{theorem}

\bigskip

The sufficiency of the theorem is a consequence of Proposition \ref{th1.3}, and the explicit expression of the above ${\mathcal D}$ follows easily from the modulo product-form. The more challenging  part is on the necessity. If $\mu= \mu_{4, {\mathcal D}}$, then the spectrum $\Lambda$ is contained in the set of zeros ${\mathcal Z}(\widehat \mu)$ of $\widehat \mu$.  By using the expression of
$
\widehat \mu (\xi)= \prod_{n=1}^\infty (b^{-1}P_{\mathcal D} (e^{-2\pi i b^{-n} \xi}))
$
where $P_{\mathcal D} (x) = \sum_{d\in {\mathcal D}} x^d$, and $\Lambda-\Lambda \subset {\mathcal Z}(\widehat \mu)$, we show by elimination that ${\mathcal D}$ must be of the form in the theorem.

\bigskip

For the organization of the paper,  we summarize and prove some  basic results concerning the geometric  and algebraic structures of the general tiles and spectral sets in Section 2. In Section 3, we introduce the various type of  product-forms, and establish the spectral property of the corresponding self-similar tiles. Theorem \ref{th1.1} is proved there. The proof of Theorem \ref{th1.2} is longer and is finished in
Section 4.  Theorem \ref{th1.4} is proved  in Section 5. We conclude with some remarks in Section 6.

\bigskip

\section{\bf Preliminaries}

\medskip

Throughout the paper, we restrict our consideration in ${\mathbb R}$ even though some of the results can be stated in ${\mathbb R}^s$. Let $\mathcal L$ denote the Lebesgue measure. Let $K$ be a bounded Borel measurable subset in ${\mathbb R}$, we use ${\mathcal L}|_K$ to denote the restriction of $\mathcal L$ on $K$, and $|K|$ the Lebesgue measure of $K$. For $A, B \subset {\mathbb R}$, $A\oplus B = \{a+b: a\in A, b\in B\}$, with $\oplus$ meaning  all elements in the sum are distinct. For $A = \{x_1, \dotsc, x_k\}$,
we use $A  (\hbox {mod} \ n)$  to mean the set $\{x_1 + nt_1, \dotsc, x_k + nt_k\}$ for some $ t_1, \dotsc, t_k \in {\mathbb Z}$.

\medskip

 We say that $K$ is a {\it tile} if there exists a countable subset ${\mathcal J} \subset {\mathbb R}$ such that $K +{\mathcal J} = {\mathbb R}$ a.e. and $\{K + t: \ t \in {\mathcal J}\}$ is a disjoint family up to sets of Lebesgue measure zero.
We call ${\mathcal J}$ a {\it tiling set} for $K$,  and  call it {\it periodic}  if there exists $0<\tau \in {\mathbb R}$  such that ${\mathcal J} = \tau + {\mathcal J}$. It is well-known that a tile $K \subset{\mathbb R}$ always admits a periodic tiling set $\mathcal J$ of the form
\begin{eqnarray*}
{\mathcal J} = B \oplus \tau\mathbb Z,
\end{eqnarray*}
where $B$ is a finite set of rational numbers, and $\tau = \#B \cdot |K|$ \
is the period of ${\mathcal J}$ [LW2, Theorems 1 and
3]. It follows from [LW2,
Theorems 3] that there exists an integer $n$ such that
${\mathcal B}: = n B/\tau \subset {\mathbb Z}_n:=\{0,1,\dotsc,
n-1\}$  and ${\mathcal B}$ tiles ${\mathbb Z}_n$, i.e.,  there
exists ${\mathcal A}\subset \mathbb Z$ such that
$$
{\mathcal A}\oplus{\mathcal B}\equiv{\mathbb Z}_n (\text{mod} \ n).
$$
Hence, by rescaling $K$  (to $|K| = n/\#{\mathcal B}= \#{\mathcal
A}$), we can always assume that ${\mathcal J}\subseteq {\mathbb Z}$
and has period $n>0$, and ${\mathcal J}$ has the form
\begin {equation} \label {eq2.1}
{\mathcal J}= {\mathcal B}  \oplus n{\mathbb Z}
\end{equation}
for some finite subset ${\mathcal B} \subset {\mathbb Z}_n$.  Indeed, this is the canonical tiling set for the self-similar tiles we consider.
A tile $K$ admits a partition as follows.   Define, for each $x \in [0,1]$, $$
{\mathcal A}_x = \{j\in {\mathbb Z} : \ x + j \in K\}.
$$
Then, as $K$ is bounded, there are only finitely many $\A_x$. Consider
$x\sim y$ if $\A_x= \A_y$, then there are only finitely many
equivalent classes $[x]$ with positive Lebesgue measure. We list the
corresponding $\A_x$  as $\A_1, \dotsc , \A_k$, and define
$$
K_j = \{x \in [0,1]: \ \A_x = \A_j\}.
$$

\bigskip

\begin{theorem}[\cite{LW2}, Theorem 3]\label{th2.1}
Suppose $K$ is a tile in ${\mathbb R}$ with a tiling set ${\mathcal J}={\mathcal B}\oplus n{\mathbb Z}$ as in (\ref{eq2.1}). Then $K$ has a partition $K=\bigcup_{j=1}^{k}(K_j+{\mathcal A}_j)$ up to a set of measure zero,  and $K_j$ and ${\mathcal A}_j \subset {\Bbb Z}$ satisfy

\vspace{0.1cm}
\ (i)  ${\mathcal A}_j\oplus {\mathcal B}\equiv{\mathbb Z}_n (\text{mod} \ n)$ for $j=1,\dotsc , k$,

\vspace{0.15cm}

(ii) $\bigcup_{j=1}^{k}K_j=[0,1]$ a.e.
\end{theorem}

\bigskip

Next, we recall some basic facts on spectral sets. Let $\mu$ be a bounded regular Borel measure on ${\mathbb R}$ with compact support $K$, and let $\widehat \mu$ denote its Fourier transform, $\widehat \mu(\xi) = \int_{\mathbb R} e^{-2 \pi i \xi x}d \mu(x)$. A basic criterion for $\mu$ to be  a spectral measure is that [JP]: {\it there exists a countable set  $\Lambda \subset \R$ (spectrum) such that
\begin{equation}\label {eq2.2}
\sum_{\lambda\in\Lambda}|\widehat\mu(\xi+\lambda)|^2 = 1 \quad  {\text{for}} \ a.e. \ \xi\in{\mathbb R}.
\end{equation}}

For a spectral set $K$  (i.e., $\mu$ is the Lebesgue measure on $K$) with $|K|=1$, we have the following structural theorem for the spectrum  $\Lambda$ ([P], [DJ2]):
there exists an integer $p>0$ such that  $\Lambda=\{0,\lambda_1,\dotsc,\lambda_{p-1}\}+p\mathbb{Z}$ with $\lambda_i\in(0,p)$.
Moreover, $K$ $p$-tiles $\mathbb{R}$ by ${\mathbb Z}/p$-translations, i.e., $K$ satisfies
$$
\sum_{j\in\mathbb{Z}}\chi_K(x+{j}/{p})=p \quad  \hbox {for} \ a.e. \  x\in[0,{1}/{p}].
$$
(It is weaker than being a tile.) The spectral set $K$ also admits a decomposition analogous to the tile case in Theorem \ref{th2.1}, which reveals the close relationship of the two. Specifically, for $x\in{\mathbb R}$,  we define
$
J_x=\{j\in {\mathbb Z}: \ x+{j}/{p}\in K\}=\{j_1,j_2\dotsc,j_p\}.
$ There are finitely many such $J_x$, we let $J_1, \dotsc , J_m$ be those $J_x$ such that
$
K_j:=\{x\in[0,{1}/{p}):\ J_x=J_j\}
$
has positive Lebesgue measure. These $K_j$'s  are disjoint up to a set of Lesbegue measure zero.

\medskip

\begin{theorem}[\cite {P,DJ2}]\label{th2.2}
Let $K$ be a spectral set with Lebesgue measure $1$. Then the spectrum $\Lambda$ has the form $\Lambda=\{0,\lambda_1,\dotsc,\lambda_{p-1}\}+p{\mathbb Z}$ with $\lambda_i\in(0,p)$. Moreover, for such $\Lambda$, $K$ has a partition up to a measure zero set of the form $
K=\bigcup_{j=1}^{m}(K_j+J_j/p),$
where $K_j$ and $J_j$ satisfy

\vspace{0.15cm}

\ (i) For $j=1,\dotsc , m$, $J_j\subset{\mathbb Z}$  and $J_j/p$ have a common spectrum $\{0,\lambda_1,\dotsc,\lambda_{p-1}\}$;

\vspace{0.15cm}

(ii)  $\bigcup_{j=1}^{m}K_j=[0,1/p]$ a.e.

\end{theorem}

\bigskip

We will need the following two propositions in the sequel. The first one was proved independently in \cite{LW4} and \cite {P}. Note that the necessity is a direct consequence of the expression of the spectrum in Theorem \ref{th2.2}.

\medskip

\begin{Pro}\label{th2.3}
Let $K \subset {\mathbb R}$ be a Borel measurable subset with finite positive Lebesgue measure. Suppose $0\in \Gamma\subset{\mathbb R}$ is a finite set with $(\Gamma-\Gamma)\cap {\mathbb Z}=\{0\}$ and $\{e^{2\pi i\lambda \cdot x}: \lambda\in{\Gamma + \mathbb Z}\}$ is orthogonal in $L^2(K)$. Then $ \Gamma + {\mathbb Z}$ is a spectrum of $K$ if and only if $\#\Gamma=|K|$.
\end{Pro}

\bigskip

If $K$ is a tile as well as a spectral set, then as a consequence of the above two theorems, we have

\medskip

\begin{Pro}\label{th2.4}
Let $K \subset\mathbb R$ be a  tile with a tiling set ${\mathcal J}= {\mathcal B} \oplus n{\mathbb Z}$ as in (\ref{eq2.1}). Assume that $K$
is also a spectral set. Then $K$ admits a spectrum of the form $\Lambda= \Gamma \oplus {\mathbb Z}$ with $\#\Gamma  = |K|$.
Moreover, if ${\mathcal A}_j$ is such that $ {\mathcal A}_j\oplus {\mathcal B}\equiv {\mathbb Z}_n \ (mod \ n)$ as in Theorem \ref{th2.1}, then all the ${\mathcal A}_j$'s have a common spectrum $\Gamma$.

\end{Pro}

\medskip

\noindent {\bf Proof.} Let $\ell =|K|$, then the tiling property implies $|K| = \#{\mathcal A}_j$ for all $j$, so that $\ell$ is a positive integer.  Let $\widetilde K= K/\ell$ and let $\widetilde \Lambda$ be the spectrum of $\widetilde K$.  By the tiling property, we have
\begin{eqnarray}\label{eq2.3}
\sum_{j\in\mathbb{Z}}\chi_{\widetilde{K}}(x+{j}/{\ell})=\ell \quad  \hbox {for} \  a.e. \ x\in[0,{1}/{\ell}].
\end{eqnarray}
For the function $\frac{1}{\ell}\sum_{j\in\mathbb{Z}}
\chi_{\widetilde{K}}(x+{j}/{\ell})$ on $[0, 1/\ell]$, its Fourier coefficients satisfy
\begin{eqnarray*}
\delta_{k,0} = a_k&=&\ell\int_{0}^{\frac{1}{\ell}}\frac{1}{\ell}
\sum_{j\in\mathbb{Z}}\chi_{\widetilde K}(x+{j}/{\ell}) e^{-2\pi i\ell kx} \ dx\\
&=&\int_{\mathbb R}\chi_{\widetilde K}(x) e^{-2\pi i\ell kx} \ dx=\widehat{\chi}_{\widetilde K}(\ell k), \quad k \in {\mathbb Z}.
\end{eqnarray*}
This shows that $\ell \mathbb Z$ is an orthonormal set in $L^2(\widetilde K, \mu)$, where $\mu = {\mathcal L}|_{\widetilde K}$, and thus $\ell\mathbb Z\subset \widetilde{\Lambda}$. From  Theorem \ref{th2.2},
\begin{eqnarray*}
\widetilde{\Lambda}=\{0,\lambda_1,\dotsc,\lambda_{p-1}\}\oplus p{\mathbb Z}, \ \text{with} \ \lambda_{i}\in(0,p).
\end{eqnarray*}
This implies $p|\ell$. We write  $\ell=pm$, and let
$$
\widetilde \Gamma = \{0,\lambda_1,\dotsc ,\lambda_{p-1}\}\oplus p\{0, 1,\dotsc, m-1\},
$$
it follows that  $\widetilde \Lambda = \widetilde \Gamma \oplus \ell {\mathbb Z}$. Let $\Gamma = \frac 1 {\ell} \widetilde\Gamma$, then $\Lambda = \Gamma \oplus {\mathbb Z}$ is a spectrum of $K$.

\medskip

Note that for $\widetilde K$ with spectrum $\widetilde \Gamma$, the family of $\A_j/\ell$ coincides with the family of $J_j/\ell$ in Theorem \ref{th2.2}, hence, the last statement follows. \qquad $\Box$

\bigskip

We conclude this section by recalling a result on the tiling and spectral property of certain finite sets. A finite subset $\A\subset \Z$ is called an {\it integer tile} if there exists an $n$ and a finite set $\B\subset \mathbb Z$ such that
\begin{equation} \label{eq2.4}
\A\oplus \B \equiv \Z_n (\hbox {mod} \ n)
\end{equation}
In this case $\A$ tiles $\Z$  (and hence $[0,1] + \A$ tiles ${\mathbb R}$)  with the tiling set $\B + n \Z$ as  in (\ref{eq2.1}).  The following theorem is due to Pedersen and Wang [PW] (see also [DJ1]). A short proof using a cyclotomic polynomial approach is in [HLL]. Note that the condition on $\A$ is stronger than being an integer tile as in (\ref{eq2.4}).

\bigskip

\begin{theorem} \label {th2.5}
 Let $\A \subset {\mathbb Z}^+$ be a finite set with $0\in \A$. Suppose there exists $n$ and $\B \subset \Z^+$ such that
$
\A\oplus \B = \Z_n.
$
Then $\A$ is a spectral set with a spectrum $\Lambda \subset \Z/n$.
\end{theorem}

\bigskip
\bigskip

\section {\bf Spectral property of product-forms}

\medskip

In the following, we assume, without loss of generality, that $0 \in {\mathcal D}\subset {\mathbb Z}^+$ and gcd$({\mathcal D})=1$. Let $b = \#{\mathcal D}$ and let $\Z_b = \{0, \dotsc , b-1\}$.  For $K= K(b, {\mathcal D})$, the invariant identity (\ref{eq1.1}) reduces to
\begin{equation} \label {eq3.1}
bK = K+ {\mathcal D}.
\end{equation}
The corresponding self-similar measure $\mu : = \mu_{b, {\mathcal D}}$ is  $\mu (\cdot ) = b^{-1} \sum_{d\in {\mathcal D}}\mu (b\cdot -d)$,  and its Fourier transform is
\begin{equation} \label {eq3.2}
\widehat\mu (\xi) = \prod_{n=1}^{\infty} \big ( b^{-1} P_{{\mathcal D}}(e^{-2\pi i \xi/b^n})\big ),
\end{equation}
 where $P_{{\mathcal D}}(x)= \sum_{d\in {\mathcal D}}x^d$. We call $P_{{\mathcal D}}$ the {\it mask polynomial} of ${\mathcal D}$.

\medskip

It is well-known that if ${\mathcal D} \equiv {\mathbb Z}_b \ (\text{mod} \ b)$, a complete residue set of $b$, then it is a tile digit set.  In this case, $K(b, {\mathcal D})$ has Lebesgue measure $1$, and has ${\mathbb Z}$ as a tiling set; also  the dual lattice ${\mathbb Z}$ is a spectrum of $K$. The following is an important extension of the complete residue sets.

\bigskip

\begin{Def}\label{th3.1}\cite{LW1}
We call ${\mathcal D}$  a product-form digit set with respect to $b$ if
\begin{eqnarray}\label{eq3.3}
{\mathcal D}={\mathcal E}_0\oplus b^{\ell_1}{\mathcal E}_1\oplus\dotsb\oplus b^{\ell_k}{\mathcal E}_k,
\end{eqnarray}
where ${\mathcal E}={\mathcal E}_0 \oplus {\mathcal E}_1\oplus \dotsb \oplus {\mathcal E}_k\equiv {\mathbb Z}_b\ (\text{mod}\ b)$, and $0\le \ell_1 \le \ell_2 \le \dotsm\le\ell_k$. If ${\mathcal E}=\Z_b$, then ${\mathcal D}$ is called a strict product-form \cite{O}.
\end{Def}

\bigskip
The basic property of the product-form digit set is contained in the following proposition.

\medskip

 \begin{Pro}\label{th3.2}
For the product-form ${\mathcal D}$ in (\ref{eq3.3}),
let
\begin{eqnarray}\label{eq3.4}
{\mathcal A} =  \bigoplus_{i=1}^{k}\bigoplus_{j=0}^{\ell_i-1} b^j{\mathcal E}_i .
\end{eqnarray}
 Then ${\mathcal A}$ is an integer tile,  $K(b, {\mathcal D})$ is a self-affine tile satisfying
\begin{eqnarray}\label{eq3.5}
K(b,{\mathcal D})=K(b,{\mathcal E})\oplus{\mathcal A}
\end{eqnarray}
and $|K(b,  {\mathcal D})| = \# {\mathcal A}$.
\end{Pro}

\medskip

\noindent {\bf Proof.} That $K(b, {\mathcal D})$ is a tile and satisfies (\ref {eq3.5}) was proved in [LW1].  It follows easily by checking $K(b,{\mathcal E})\oplus{\mathcal A}$ satisfies the self-similar identify (\ref{eq3.1}).
To show that $\A$ is an integer tile, we let $\ell_0=0$ and
\begin{eqnarray}\label{eq3.6}
{\mathcal B} =  \bigoplus_{i=0}^{k-1}\bigoplus_{j=\ell_i}^{\ell_k-1}b^j{\mathcal E}_i.
\end{eqnarray}
 Then
\begin{eqnarray} \label {eq3.7}
{\mathcal A} + \ {\mathcal B}
& = & \big(\bigoplus_{j=0}^{\ell_{k}-1} b^j \E_0 \big)\oplus \dotsc \oplus \big(\bigoplus_{j=0}^{\ell_{k}-1} b^j \E_k \big)  \\
& =& \ {\mathcal E} \oplus b{\mathcal E} \oplus \dotsc \oplus b^{\ell_k-1}{\mathcal E} \ \equiv \ {\mathbb Z}_{b^{\ell_k}} \ (\text{mod} \ b^{\ell_k}). \nonumber
\end{eqnarray}
 The above sum is actually direct since
$\#\A \cdot \#\B = \prod_{i=1}^k \#\E_i^{\ell_i} \prod_{i=0}^{k-1} \#\E_i^{\ell_k -\ell_i} = b^{\ell_k}$. This shows that ${\mathcal A}$ tiles ${\mathbb Z}_{b^{\ell_k}}$, i.e. ${\mathcal A}$ is an integer tile.
 On the other hand, the assumption that $\text{gcd}({\mathcal D})=1$ implies  $\text{gcd}({\mathcal E})=1$. Hence, $|K(b,{\mathcal E})|=1$
follows from ${\mathcal E}\equiv {\mathbb Z}_b (\text{mod} \ b)$. Therefore, we have $|K(b,  {\mathcal D})| = \# {\mathcal A}$.
 \qquad $\Box$

\bigskip

 Note that an integer tile is weaker than being a tile digit set and the above  ${\mathcal A}$ is not necessarily a tile digit set.
For example, let ${\mathcal D}={\mathcal E}_0\oplus 8{\mathcal E}_1$ with ${\mathcal E}_0=\{0, 2\}$ and ${\mathcal E}_1=\{0, 1, 4, 5\}$. Then ${\mathcal D}$ is a strict product-form, and ${\mathcal A}={\mathcal E}_1=\{0, 1, 4, 5\}$. It is easy to show that ${\mathcal A}$ is not a tile digit set with respect to $4$ since $|K(4, {\mathcal A})|=0$ by checking the identity (\ref{eq3.1}).

 \medskip

\begin {theorem} \label {th3.3}
Suppose ${\mathcal D}$ is a strict product-form digit set, then
$K(b,{\mathcal D})$ is a finite union of intervals $[0,1]$, and is a spectral set with spectrum $\Gamma \oplus {\mathbb Z} \subseteq {\mathbb Z}/{n}$ for some $n>0$.
\end{theorem}

\medskip

 \noindent {\bf Proof.}\ If ${\mathcal D}$ is a strict product-form, then $\E =\Z_b$. Hence  $K(b,{\mathcal E})=[0,1]$  and the corresponding self-affine tile satisfies $K(b,{\mathcal D})=[0,1]\oplus{\mathcal A}$. This implies that $K := K(b, {\mathcal D})$ is a finite union of intervals. Let ${\mathcal B}$ be as in (\ref{eq3.6}), then by $\E =\Z_b$ and a direct calculation as in (\ref{eq3.7}),  we obtain
$$
{\mathcal A}\oplus {\mathcal B}=\{0, 1,\dotsc, b^{\ell_k}-1\}.
$$
 Theorem \ref{th2.5} implies  that $\A$  is a spectral set with a spectrum $\Gamma\subset {\mathbb Z}/{n}$.  Hence $\Gamma \oplus {\mathbb Z} \subseteq {\mathbb Z}/{n}$ gives an exponential orthonormal set of $L^2(K)$.  As $|K| = \#\A$, it follows from Proposition \ref{th2.3} that  $\Gamma \oplus {\mathbb Z} $ is actually the spectrum of $K$.
\qquad $\Box$

\bigskip

\begin{theorem}\label{th3.4}
Let ${\mathcal D}$ be a product-form digit set. Then $ K(b, {\mathcal D})$ is a spectral set if and only if ${\mathcal A}$ is a spectral set.
\end{theorem}

\medskip

\noindent {\bf Proof.} \
Since ${\mathcal E}$ is a complete residue set modulo $b$, $K(b, {\mathcal E})$ has ${\mathbb Z}$ as a tiling set, and the spectrum is ${\mathbb Z}$.  Let $\delta_\A$ be the point mass measure of $\A$, it follows from  (\ref{eq3.5}) that
\begin{eqnarray*}
\mu_{b,{\mathcal D}}={\mathcal L}|_{K(b,{\mathcal D})}
={\mathcal L}|_{K(b,{\mathcal E})}\ast\delta_{{\mathcal A}}=\mu_{b,{\mathcal E}}\ast\delta_{{\mathcal A}}.
\end{eqnarray*}
Assume that $ K = K(b, {\mathcal D})$ is a spectral set with spectrum $\Lambda$. By Proposition \ref{th2.4}, we have
$\Lambda=\Gamma \oplus \mathbb Z$ with $(\Gamma-\Gamma)\cap {\mathbb Z}=\{0\}$ and $\#\Gamma=|K|$. Then for $\Lambda=: \Gamma \oplus {\mathbb Z}$,
\begin{eqnarray}\label{eq3.8}
1=\sum_{\lambda\in\Lambda}|\widehat\mu_{b,{\mathcal D}}(x+\lambda)|^2
&=&\sum_{\lambda\in\Lambda}|\widehat\mu_{b,{\mathcal E}}
(x+\lambda)|^2\lvert\widehat\delta_{{\mathcal A}}(x+\lambda)|^2 \nonumber \\
&=&\sum_{{k\in{\mathbb Z}},\ {\gamma\in \Gamma}}
|\widehat\mu_{b,{\mathcal E}}(x+k+\gamma)|^2 |\widehat\delta_{{\mathcal A}}
(x+k+\gamma)|^2 \nonumber \\
&=&\sum_{\gamma\in\Gamma}|\widehat\delta_{{\mathcal A}}(x+\gamma)|^2\sum_{k\in {\mathbb Z}}|
\widehat\mu_{b, {\mathcal E}}(x+k+\gamma)|^2 \nonumber\\
&=&\sum_{\gamma\in\Gamma } |\widehat\delta_{{\mathcal A}}
(x+\gamma)|^2.
\end{eqnarray}
It follows from the criterion in (\ref {eq2.2}) that $\Gamma$ is a spectrum of ${\mathcal A}$.

\vspace{0.2cm}

 Conversely, assume that ${\mathcal A}$ is a spectral set and $\Gamma$ is its spectrum. Then for $\Lambda=: \Gamma \oplus {\mathbb Z}$, using the same argument as in (\ref{eq3.8}), we have
$\sum_{\lambda\in\Lambda}|\widehat\mu_{b,{\mathcal D}}(x+\lambda)|^2=\sum_{\gamma\in\Gamma } |\widehat\delta_{{\mathcal A}}(x+\gamma)|^2=1$.
This shows that $\Lambda$ is a spectrum of $K(b, {\mathcal D})$ by the criterion in (\ref {eq2.2}). \qquad $\Box$

\bigskip

Next, we will consider the modulo-product form. For a digit set ${\mathcal D}$,  we will make use of the mask polynomial
$
P_{{\mathcal D}}(x) = \sum_{d\in {\mathcal D}} x^d.
$
It is clear that $\widehat \delta_{{\mathcal D}} (\xi) = P_{\mathcal D}(e^{-2\pi i \xi x})$.
We will deal with the direct summands in ${\mathcal D}$ by factorization of the mask polynomial through its roots of unity.  Recall that a {\it cyclotomic polynomial} $\Phi_d(x)$ is the minimal polynomial of the $d$-th root of unity, of which the explicit form can be obtained inductively by the identity
\begin{eqnarray*}
x^s-1={\prod}_{d\lvert s}\Phi_d(x).
\end{eqnarray*}
 It is clear that $\Phi_d(x)\lvert P_{\mathcal{D}}(x)$ if and only if $P_{\mathcal{D}}(e^{\frac{2\pi i}{d}})=0$.
It follows that if ${\mathcal E} = {\mathcal E}_0 \oplus \cdots \oplus {\mathcal E}_k\equiv {\mathbb Z}_b (\hbox {mod} \ b)$, then
\begin{eqnarray}\label{eq3.9}
P_{{\mathcal E}}(x)=P_{{\mathcal E}_0}(x)\dotsb P_{{\mathcal E}_k}(x)\equiv 1+x+\dotsb+x^{b-1} (\text{mod} \ (x^b-1)).
\end{eqnarray}
If ${\mathcal D}$ is a product-form in (\ref{eq3.3}), then its mask polynomial can be expressed as
\begin{eqnarray}\label{eq3.10}
P_{\mathcal D} (x)=P_{{\mathcal E}_0}(x)P_{{\mathcal E}_1}(x^{b^{\ell_1}}) \dotsb P_{{\mathcal E}_k}(x^{b^{\ell_k}}).
\end{eqnarray}
We let  $S_i=\{d>1: \ d| b,\ \Phi_d(x)| P_{{\mathcal E}_i}(x)\}$,  and let
$\Psi_i(x)=\prod_{d\in S_i}\Phi_d(x)$.
Then $\Psi_i(x)| P_{{\mathcal E}_i}(x)$, and $\Psi_i(x^{b^{\ell_i}})\mid P_{{\mathcal E}_i}(x^{b^{\ell_i}})$. Let
\begin{eqnarray}\label{eq3.11}
K_{\mathcal{D}}^{(i)}(x)=\Psi_0(x)\Psi_1(x^{b^{\ell_1}})
\dotsb\Psi_i(x^{b^{\ell_i}}), \ \ 0\le i\le k,
\end{eqnarray}
and let $K_{{\mathcal D}}(x)=K_{{\mathcal D}}^{(k)}(x)$. Then $K_{{\mathcal D}}(x)\mid  P_{\mathcal{D}}(x)$.

\bigskip

\begin{Def}[\cite{LLR1}]\label{th3.5}
For ${\mathcal E} = {\mathcal E}_0 \oplus \cdots \oplus {\mathcal E}_k\equiv {\mathbb Z}_b \ (\hbox {mod} \ b)$, $0\leq \ell_1\leq \cdots \leq \ell_k$,   let $
n_i=\text{l.c.m.}\{s:\Phi_s(x)\mid K_{\mathcal{D}}^{(i)}(x)\}$. Define ${\mathcal D}^{(0)}\equiv {\mathcal E}_0 \ (\text{mod} \ n_0)$, and
$$
{\mathcal D}^{(i)}\equiv {\mathcal D}^{(i-1)}\oplus b^{\ell_i}{\mathcal E}_i \ (\text{mod} \ n_i), \quad  1\leq i \leq k\ .
$$
We say that $\mathcal{D}$ is a modulo product-form if ${\mathcal D}={\mathcal D}^{(k)}$.
\end{Def}

\bigskip
 That a modulo product-form is a tile digit set is proved in [LLR1, Theorem 3.5]. It is clear that a product-form is a modulo product-form by ignoring all the modulo actions. Our main result is

\bigskip

\begin{theorem}\label{th3.6}
Let ${\mathcal D}$ be a modulo product-form with respect to $b$, and be defined by ${\mathcal E} ={\mathcal E}_0 \oplus \cdots \oplus  {\mathcal E}_k={\mathbb Z}_b$ and $0\le\ell_1\leq \cdots \leq \ell_k$. Then the self-similar tile $K(b,{\mathcal D})$ is a spectral set, and any spectrum of  the product form $K(b,{\mathcal D}^{\prime})$ with ${\mathcal D}^{\prime}={\mathcal E}_0 \oplus \cdots \oplus b^{\ell_k}{\mathcal E}_k$ is also a spectrum of $K(b,{\mathcal D})$.
\end{theorem}

\medskip

Since the proof of theorem involves a few technical steps, we will carry it out separately in the next section. In the following, we will apply the theorem to two special cases.

\medskip

For  $b= p^\alpha$  or $b= pq$ where $p, q$ are distinct primes, it is known that  ${\mathcal D}$ is a tile digit set if and only it is a modulo product-form (\cite{LLR2}, see also \cite{LW1, LR}).

\bigskip

\begin {Pro} \label {th3.7}
Let ${\mathcal D}$ be a digit set with $\#{\mathcal D} = b = p^\alpha$ or $pq$, where $p, q$ are distinct primes. Then  $K(b, {\mathcal D})$ is a spectral set.
\end{Pro}

\medskip

\noindent {Proof.} We first consider  $b= p^\alpha$.  Let ${\mathcal E}_i= p^{i} \{0, \dotsc, p-1\}$, and note that ${\mathcal E} = {\mathcal E}_0 \oplus \cdots \oplus {\mathcal E}_{\alpha -1} = {\mathbb Z}_{p^\alpha}$.
 Then, by [LLR2, Theorems 4.3, 5.1], the tile digit set is a modulo product-form with respect to ${\mathcal E}_0 \oplus \cdots \oplus {\mathcal E}_{\alpha -1}$.
It follows from
Theorem \ref{th3.6} that $K(p^{\alpha}, {\mathcal D})$ is a spectral set.

\medskip
For $b=pq$, the tile digit sets are characterized by
$$
{\mathcal D} = \big ( {\mathcal E}_1 (\hbox {mod} \ p) \oplus b^{n-1} {\mathcal E}_2\big )\ (\hbox {mod} \  b^n)
$$
where ${\mathcal E}_1 = \{0, \dotsc, p-1\}$ and ${\mathcal E}_2 = \{0, p, \dotsc, p(q-1)\}$ (or interchange the role of $p,q$) \cite{LR}. Note that ${\mathcal E}_1 \oplus {\mathcal E}_2= {\mathbb Z}_{pq}$, and the corresponding product-forms  are strict. Similar to the above, we conclude that such $K(b, {\mathcal D})$ are spectral sets. \qquad $\Box$

\bigskip

As an illustration, we consider the case $\#{\mathcal D} =4$.

\medskip

\begin {Example} \label{th3.8}
Let  ${\mathcal D} \subset{\mathbb Z}^+$ be a digit set with $\#{\mathcal D} = 4$ and $\text{gcd} ({\mathcal D})=1$. Then ${\mathcal D}$ is a tile digit set with respect to $4$ if and only if ${\mathcal D} = \{0, a, 2^{t} \ell,  a +2^{t} \ell'\}$, where  $a, t, \ell, \ell'$ are odd integers. In this case, $K(4, {\mathcal D})$ is a spectral set.
\end{Example}

\medskip

\noindent {\bf Proof.}  By the characterization of the tile digit set $\mathcal D$ with $\#{\mathcal D}=p^\alpha = 2^2$  as the modulo product-form,  ${\mathcal D}$ has the form
$$
{\mathcal D} = \big (\{0,1\}  (\hbox {mod}\ 2) \oplus 2^{2t_1} \{0,2\}\big ) (\hbox {mod}\ 2^{2(t_1+1)}).
$$
The explicit expression of ${\mathcal D}$ in the theorem can be derived easily from this form. It is a spectral set by Proposition \ref {th3.7}.   \qquad $\Box$

\bigskip

We remark that the converse, i.e.,  $K(4, {\mathcal D})$ is a spectral set implies that ${\mathcal D} = \{0, a, 2^{t} \ell,  a +2^{t} \ell'\}, \ a, t, \ell, \ell'$ are odd integers, is more complicate. We will proof  this in detail in Section 5.

\bigskip

\section {\bf Proof of Theorem \ref{th3.6}}

\medskip

Let ${\mathcal D}$ be a modulo product-form defined by ${\mathcal E}={\mathcal E}_0\oplus \dotsm \oplus {\mathcal E}_k={\mathbb Z}_b$ and $0\le \ell_1\le \cdots \le \ell_k$ as in Definition \ref{th3.5},  and let ${\mathcal D}^\prime$ be the associated strict product-form.
The strategy of proving Theorem \ref{th3.6} is to show $|K(b,{\mathcal D})|=|K(b,{\mathcal D}^{\prime})|$ first, then use Propositions \ref{th2.3} and \ref{th2.4} to conclude the statement. To prove the identity, we make use of the fact that the Lebesgue measure of an integral tile $K\subset {\mathbb R}$  can be obtained from its periodic tiling set (see Section 2). For this, we recall a lemma concerning the tiling set of a self-similar tile in ${\mathbb R}$ (Theorem 3.1(ii) in [LR]).

\begin{Lem}\label{th4.1}
Let $0\in {\mathcal D}\subset {\mathbb Z}$  and $\text{gcd}({\mathcal D})=1$. Suppose $\D$ is a tile digit set and ${\mathcal J}\subseteq {\mathbb Z}$ is periodic and
${\mathcal J}=b{\mathcal J}\oplus {\mathcal D}$. Then ${\mathcal J}$ is a tiling set for $K(b, {\mathcal D})$. In this case, ${\mathcal J}$ is called {\it a self-replicating tiling set}.
\end{Lem}

Suppose that $b= p_0\dots p_m$, where $p_i, \ 0\le i\le m,$ are primes. Let $\widetilde p_i=  p_0\dotsm p_i$, and let $\C_0={\mathbb Z}_{p_0}, \C_i=\tilde p_{i-1}{\mathbb Z}_{p_i}$ for $1\le i\le m$. Then, by applying de Brujin's decomposition theorem [dB] inductively (see Proposition 2.1 and Corollary 2.2 in \cite{PW}, starting with ${\mathbb Z}_b = \{0\} \oplus {\mathbb Z}_b$), we have

\begin{eqnarray}\label{eq4.1}
{\mathbb Z}_b=\C_0\oplus \C_1 \oplus \dotsm \oplus \C_m.
\end{eqnarray}
Note that each component is ``irreducible", and  all irreducible decomposition of ${\mathbb Z}_b$ must be of the above form with respect to some rearrangement $b= p'_0\cdots p'_m$.

\medskip

\begin{Lem}\label{th4.2}
For any integer $b\geq 2$,
if \ ${\mathbb Z}_b$ is decomposed as ${\mathbb Z}_b={\mathcal E}_0 \oplus \cdots \oplus {\mathcal E}_k$, then  ${\mathcal E}_i$ has the  form:
$$
{\mathcal E}_i=\C_{i_1}\oplus \C_{i_2}\oplus \dotsm\oplus \C_{i_{s_i}} \ \text{with} \  0\le i_1< i_2<\dotsm<i_{s_i}\le m,
$$
where $\C_0={\mathbb Z}_{p_0}$ and $\C_i=\widetilde p_{i-1}{\mathbb Z}_{p_i}$, $1\le i\le m$, for some arrangement $b = p_0 \cdots p_m$.
\end{Lem}

\bigskip

\noindent {\bf Proof.} It suffices to consider  ${\mathbb Z}_b={\mathcal E}_0 \oplus {\mathcal E}_1$. Let $b= p'_0\cdots p'_m$. As ${\mathbb Z}_b$ is composed of summands $\alpha_i {\mathbb Z}_{p'_i}$ and they are irreducible, they must be summands of either $\E_0$ or $\E_1$. Since the decomposition of $\Z_b$ must be of the form in (\ref {eq4.1}), there exists an arrangement $p_0\cdots p_m=b$ with each $ \C_j, \ 0 \leq j \leq m $, equals one and only one $ \alpha_i {\mathbb Z}_{p'_i}$.   This proves the lemma.  \qquad $\Box$

\bigskip

Now let ${\mathcal E}={\mathcal E}_0 \oplus \cdots \oplus {\mathcal E}_k={\mathbb Z}_b$.  For $0= \ell_0 <\ell_1<\cdots <\ell_k$, let ${\mathcal D}^\prime$ be the product-form:
$$
{\mathcal D}^\prime = {\mathcal E}_0 \oplus b^{\ell_1}{\mathcal E}_1\oplus\cdots \oplus b^{\ell_k}{\mathcal E}_k.
$$
It is seen from the  proof in Proposition \ref {th3.2} that
$$
{\mathcal J}^\prime =b^{\ell_k}{\mathbb Z}\oplus \bigoplus_{i=0}^{k-1}\bigoplus_{j=\ell_{i}}^{\ell_k-1} b^j {\mathcal E}_i
 $$
is a tiling set for $K(b,{\mathcal D}^\prime)$, and is also self-replicating. In the notation of
Definition \ref{th3.5}, ${\mathcal D}^{\prime(i)} = \bigoplus\limits_{j=0}^i b^{\ell_j}\E_j$, we can rewrite ${\mathcal J}^\prime$  as
\begin {equation} \label {eq4.2}
{\mathcal J}^\prime =b^{\ell_k}{\mathbb Z}\oplus
\bigoplus_{i=0}^{k-1}\bigoplus_{j=\ell_k-\ell_{k-i}}^{\ell_k-\ell_{k-i-1}-1}b^j {\mathcal D}^{\prime(k-i-1)}.
\end{equation}

Let $\mathcal D$ be the corresponding modulo product-form, it follows from the definition of $n_i$ in Definition \ref{th3.5} that
$ n_i=b^{\ell_i}\widetilde p_{i_{s_i}}.$
Note that the value of $n_i$ is determined by  the last summand $b^{\ell_i}\C_{i_{s_i}}$ of $b^{\ell_i}{\mathcal E}_i$.
We will show that the expression in (\ref{eq4.3}) also provides a  self-replicating tiling set ${\mathcal J}$ for $K(b, {\mathcal D})$.
We consider the tiling property of $K(b,{\mathcal D})$ in the following two lemmas: the case $n_k = b^{\ell_k+1}$ and the case $n_k < b^{\ell_k+1}$. To simplify the notations, we let $t_0 = 0$,
$
t_i = \ell_k - \ell_{k-i}$.
\bigskip

\begin{Lem}\label{th4.3}
If  $n_k=b^{\ell_k+1}$, then
\begin{equation}\label {eq4.3}
{\mathcal J} = b^{\ell_k}{\mathbb Z}\oplus
\bigoplus_{i=0}^{k-1} \bigoplus_{j=t_i}^{t_{i+1}-1} b^j{\mathcal D}^{(k-i-1)}.
\end{equation}
is a self-replicating tiling set for $K(b,{\mathcal D})$.
\end{Lem}

\medskip

\noindent {\bf Proof.}
 In view of  ${\mathcal D} := {\mathcal D}^{(k-1)}\oplus b^{\ell_k}{\mathcal E}_k \ (\hbox {mod}\ b^{\ell_k+1})$,  we write
\begin{eqnarray}\label{eq4.4}
 b{\mathcal J}\oplus ({\mathcal D}^{(k-1)}\oplus b^{\ell_k}{\mathcal E}_k) = {\mathcal J}_0 \oplus \bigoplus_{i=0}^{k-1} \bigoplus_{j=t_i +1}^{t_{i+1}-1 }b^j{\mathcal D}^{(k-i-1)} \oplus {\mathcal D}^{(k-1)}
\end{eqnarray}
where
$$
{\mathcal J}_0 = b^{\ell_k+1}{\mathbb Z}\oplus\big (\  b^{\ell_k}{\mathcal E}_k\ \oplus \
\bigoplus_{i=0}^{k-1} b^{t_{i+1}} {\mathcal D}^{(k-i-1)}\big ).
$$
For each $0\le i\le k-1$, ${\mathcal D}^{(i)}=({\mathcal D}^{(i-1)}\oplus b^{\ell_i}{\mathcal E}_i) (\text{mod} \ n_i)$. We will prove that  ${\mathcal J}_0$ has a period $b^{\ell_k}$, and  remove the modulo operation of each ${\mathcal D}^{(i)}$ step by step.

\medskip

Since $n_k=b^{\ell_k+1},$ then $\C_m$ is a summand of ${\mathcal E}_k$ and $\C_{k_{s_k}}=\C_m$. Note that
\begin{eqnarray}\label{eq4.5}
b^{\ell_k+1}{\mathbb Z}\oplus \C_m = b^{\ell_k+1}{\mathbb Z}\oplus \widetilde p_{m-1}{\mathbb Z}_{p_m}
=b^{\ell_k}\widetilde p_{m-1}(p_m{\mathbb Z}\oplus{\mathbb Z}_{p_m})
=b^{\ell_k}\widetilde p_{m-1}{\mathbb Z}.
\end{eqnarray}
 This implies that $b^{\ell_k+1}{\mathbb Z}\oplus b^{\ell_k}{\mathcal E}_k$ (and hence ${\mathcal J}_0$) has a period $b^{\ell_k}\widetilde p_{m-1}$.

 \medskip

Next, we consider $\C_{m-1}$, it is  a summand of either (i) ${\mathcal E}_k$, or (ii) ${\mathcal E}_i$, $0\leq i \leq k-1$.
In case (i), $ \C_{m-1}= \C_{k_{s_k}-1}$. By observing that ${\mathcal E}_k = \bigoplus\limits_{j=1}^{s_k}\C_{k_j}$ and
\begin {eqnarray} \label {eq4.6}
b{\mathbb Z}\oplus  \C_{m-1}\oplus \C_m
=\widetilde p_{m-2}(p_{m-1}p_m{\mathbb Z}\oplus {\mathbb Z}_{p_{m-1}}\oplus p_{m-1}{\mathbb Z}_{p_m}) = \widetilde p_{m-2}{\mathbb Z},
\end{eqnarray}
we have
\begin {eqnarray} \label{eq4.7}
b^{\ell_k+1}{\mathbb Z}\oplus b^{\ell_k}{\mathcal E}_k
 = b^{\ell_k}\widetilde p_{m-2}{\mathbb Z}\oplus b^{\ell_k}\bigoplus_{j=1}^{s_k-2}\C_{k_j}.
\end{eqnarray}
It follows from (\ref{eq4.7}) that $b^{\ell_k+1}{\mathbb Z}\oplus b^{\ell_k}{\mathcal E}_k$, and hence ${\mathcal J}_0$, have a period $b^{\ell_k}\widetilde p_{m-2}$. In case (ii), $\C_{m-1}=\C_{i_{s_i}}$. Note that $n_i=b^{\ell_i}\widetilde p_{m-1}$ and using the definition of ${\mathcal D}^{(i)}$, we have
\begin{eqnarray}\label{eq4.8}
&&b^{\ell_k+1}{\mathbb Z}\oplus b^{\ell_k}{\mathcal E}_k \oplus b^{\ell_k -\ell_i}{\mathcal D}^{(i)}\nonumber\\
&=& b^{\ell_k+1}{\mathbb Z}\oplus b^{\ell_k}{\mathcal E}_k \oplus b^{\ell_k-\ell_i}\Big (({\mathcal D}^{(i-1)}\oplus b^{\ell_i}{\mathcal E}_i) (\text{mod} \ b^{\ell_i}\widetilde p_{m-1})\Big ) \nonumber\\
&=& b^{\ell_k+1}{\mathbb Z}\oplus b^{\ell_k}{\mathcal E}_k \oplus  b^{\ell_k-\ell_i}{\mathcal D}^{(i-1)} \oplus b^{\ell_k}{\mathcal E}_i \ (\text{mod} \ b^{\ell_k}\widetilde p_{m-1}) \nonumber\\
&=&  b^{\ell_k+1}{\mathbb Z}\oplus b^{\ell_k}{\mathcal E}_k \oplus b^{\ell_k-\ell_i}{\mathcal D}^{(i-1)} \oplus b^{\ell_k}{\mathcal E}_i
\end{eqnarray}
(the last equality is obtained since $b^{\ell_k+1}{\mathbb Z}\oplus b^{\ell_k}{\mathcal E}_k$ has a period $b^{\ell_k}\widetilde p_{m-1}$ by (\ref{eq4.5})). It follows from a similar argument to (\ref{eq4.7}) that
\begin{eqnarray*}
b^{\ell_k+1}{\mathbb Z}\oplus b^{\ell_k}{\mathcal E}_k\oplus b^{\ell_k}{\mathcal E}_i
=b^{\ell_k}\widetilde p_{m-2}{\mathbb Z}\oplus b^{\ell_k}\bigoplus\limits_{j=1}^{s_k-1}\C_{k_j}\oplus b^{\ell_k}\bigoplus\limits_{j=1}^{s_i-1}\C_{k_j}.
\end{eqnarray*}
Thus,
$b^{\ell_k+1}{\mathbb Z}\oplus b^{\ell_k}{\mathcal E}_k \oplus b^{\ell_k-\ell_i}{\mathcal D}^{(i)}$, and hence ${\mathcal J}_0$,   have a period $b^{\ell_k}\widetilde p_{m-2}$.

\medskip

We continue to consider $\C_{m-3}$, it can be a summand in $\E_k$, the above $\E_i$, or a third $\E_{i\prime}$. We can use the same argument to show, say in the third case, that $b^{\ell_k+1}{\mathbb Z}\oplus b^{\ell_k}{\mathcal E}_k \oplus b^{\ell_k-\ell_{i}}{\mathcal D}^{(i)} \oplus b^{\ell_k-\ell_{i^\prime}}{\mathcal D}^{(i^\prime)}$, and hence ${\mathcal J}_0$,   have a period $b^{\ell_k}\widetilde p_{m-3}$.
By going through all the $\C_j$'s, we conclude that ${\mathcal J}_0$ has a period $b^{\ell_k}$. Also it follows from (\ref{eq4.8}) that the periodicity can nullify the modulo operation, and it yields
\begin{eqnarray}\label{eq4.9}
{\mathcal J}_0
&=&b^{\ell_k+1}{\mathbb Z}\oplus b^{\ell_k}{\mathcal E}_k\oplus
\bigoplus_{i=0}^{k-2} b^{t_{i+1}} \big ({\mathcal D}^{(k-i-2)}\oplus b^{\ell_{k-i-1}}\E_{k-i-1}\big )  \oplus b^{\ell_k}{\mathcal E}_0 \nonumber\\
&=&b^{\ell_k}{\mathbb Z}\oplus \bigoplus_{i=0}^{k-2} b^{t_{i+1}} {\mathcal D}^{(k-i-2)}.
\end{eqnarray}
Substituting (\ref{eq4.9}) into the RHS of (\ref{eq4.4})  and using the fact that ${\mathcal J}$ has a period $b^{\ell_k}$, we get
$$
b{\mathcal J}\oplus {\mathcal D} = b{\mathcal J}\oplus {\mathcal D}^{(k-1)}\oplus b^{\ell_k}{\mathcal E}_k={\mathcal J}.
$$
Hence, $\mathcal J$ is a self-replicating tiling set for $K(b, {\mathcal D})$.  Therefore  Lemma \ref{th4.1} yields the statement of the lemma.
\qquad $\Box$

\bigskip

If $n_k<b^{\ell_k+1}$, then $\E_k = \C_{k_1} \oplus \cdots \oplus \C_{k_{s_k}}$,  and $\C_{k_{s_k}+1} , \cdots , \C_m$ belongs to some $\E_i$'s, $i\not = k$. Let $G$ be the direct sum of the ${\mathcal E}_i$ which contains at least one summand from the sets $\C_{k_{s_k}+1}, \dotsc, \C_m$. It is clear that $n_k = b^{\ell_k} \widetilde p_{k_{s_k}}$, and by the same proof to (\ref{eq4.6}), we have
\begin{eqnarray*}
b{\mathbb Z}\oplus(\C_{k_{s_k}+1}\oplus \cdots \oplus \C_m)=\widetilde p_{k_{s_k}}{\mathbb Z}.
\end{eqnarray*}
It follows that $b^{\ell_k}{\mathbb Z}\oplus b^{\ell_k-1}G$  has a period $b^{\ell_k-1}\widetilde p_{k_{s_k}}$.

\medskip

We will make use of this fact to modify the tiling set ${\mathcal J}$  and ${\mathcal J}_0$ in Lemma \ref {th4.3}. The basic idea is for such $\E_i$ and $\D^{(i)}( =: \D^{(i-1)}\oplus b^{\ell_i}\E_i (\hbox {mod} \ n_i))$, we shift the $\E_i$ as a summand of $G$, and absorb the $\D^{(i-1)}$ with the sum of the other $\D^{(i-1)}$'s. Specifically, we let $t_0'=0$ and
$t_i'= t_i-1$ or $=t_i$ according to ${\mathcal E}_i, \ i =0, \dotsc, k-1 $ , is a summand in $G$ or not.

\bigskip

\begin{Lem}\label{th4.4}
For the modulo product-from ${\mathcal D}$, if  $n_k<b^{\ell_k+1}$,
we let
\begin{equation} \label {eq4.10}
{\mathcal J}=b^{\ell_k}{\mathbb Z}\oplus b^{\ell_k-1}G\oplus \bigoplus_{i=0}^{k-1}\bigoplus_{j=t'_i}^{t'_{i+1}-1} b^j \D^{(k-i-1)}.
\end{equation}
Then ${\mathcal J}$ is a self-replicating tiling set for $K(b,{\mathcal D})$.
\end{Lem}

\noindent {\bf Proof.}
We proceed in the same way as Lemma
\ref {th4.3} and consider
\begin{eqnarray*}
 b{\mathcal J}\oplus ({\mathcal D}^{(k-1)}\oplus b^{\ell_k}{\mathcal E}_k)
 = {\mathcal J}_0 \oplus \bigoplus_{i=0}^{k-1} \bigoplus_{j=t'_i +1}^{t'_{i+1}-1 }b^j{\mathcal D}^{(k-i-1)} \oplus {\mathcal D}^{(k-1)}
\end{eqnarray*}
where
$$
{\mathcal J}_0 = b^{\ell_k+1}{\mathbb Z}\ \oplus \big (\  b^{\ell_k}({\mathcal E}_k \oplus G)\ \oplus \
\bigoplus_{i=0}^{k-1} b^{t'_{i+1}} {\mathcal D}^{(k-i-1)}\big ).
$$
Since $\C_{k_{s_k}}$ is a summand of ${\mathcal E}_k$, $ b^{\ell_k+1}{\mathbb Z}\oplus b^{\ell_k}({\mathcal E}_k\oplus G)$ will contain $b^{\ell_k+1}{\mathbb Z}\oplus b^{\ell_k}(\C_{k_{s_k}}\oplus \cdots \oplus \C_m)$ as a summand.
By the same argument as in (\ref{eq4.6}), we obtain
$$
b^{\ell_k+1} {\mathbb Z}  \oplus b^{\ell_k} (\C_{k_{s_k}} \oplus \cdots \oplus \C_m) = b^{\ell_k} \tilde p_{k_{s_k}-1}{\mathbb Z}.
$$
This shows that $b^{\ell_k+1}{\mathbb Z} \oplus b^{\ell_k}({\E}_k\oplus G)$ has a period $b^{\ell_k}\widetilde p_{k_{s_k}-1}$ (compare this with the conclusion following (\ref {eq4.5})).

We then follow the same argument in Lemma \ref{th4.3} to consider
$\C_{k_{s_k}-1}, \cdots , \C_0$ inductively, and  remove the modulo operation of the associated ${\mathcal D}^{(i)}$ in ${\mathcal J}_0$ in each step. We conclude that ${\mathcal J}$  has a period $b^{{\ell_k}-1}\widetilde p_{k_{s_k}} = n_k/b$,  and
$$
b{\mathcal J}+ \D = b{\mathcal J}\oplus ( {\mathcal D}^{(k-1)}\oplus b^{\ell_k}{\mathcal E}_k) ={\mathcal J}.
$$
Therefore, ${\mathcal J}$ is a self-replicating tiling set for $K(b, \D)$.
\qquad $\Box$

\bigskip

\begin{Lem}\label{th4.5}
Let ${\mathcal D}$ be a modulo product-form and ${\mathcal D}^{\prime}$ be a  product-form as before. Then
$
|K(b,{\mathcal D})| =|K(b,{\mathcal D}^{\prime})|.
$
\end{Lem}

\medskip

\noindent{\bf Proof.} We first consider the case of $\D$ in Lemma \ref {th4.3}. For the tiling set of $K(b,{\mathcal D})$ and $K(b,{\mathcal D}^{\prime})$, we refer to (\ref{eq4.2}) and (\ref{eq4.3}) respectively. Let
$$
{\mathcal B} = \bigoplus_{i=0}^{k-1} \bigoplus_{j=t_i}^{t_{i+1}-1} b^j{\mathcal D}^{(k-i-1)}\quad  \hbox {and} \quad  {\mathcal B}' = \bigoplus_{i=0}^{k-1} \bigoplus_{j=t_i}^{t_{i+1}-1} b^j{\mathcal D}^{\prime(k-i-1)}
$$
denote the respective sums. Then ${\mathcal J} = b^{\ell_k}{\mathbb Z}\oplus {\mathcal B}$ and  ${\mathcal J}' = b^{\ell_k}{\mathbb Z}\oplus {\mathcal B}'$. Note that $\#{\mathcal B}=\#{\mathcal B}'$ by the direct sum property and $\#\D^{(i)} = \#\D^{\prime (i)}$. It follows from [LW2] (see also Section 2) that
$$
|K(b, {\mathcal D})|\ =\ {b^{\ell_k}}/{\#{\mathcal B}}\ = \ {b^{\ell_k}}/{\#{\mathcal B}'} \ = \ |K(b, {\mathcal D}')|.
$$

For the case of $\D$ in Lemma \ref {th4.4}, we define
$$
{\mathcal B} = b^{\ell_k-1}G\oplus \bigoplus_{i=0}^{k-1}\bigoplus_{j=t_i'}^{t_{i+1}'-1} b^j \D^{(k-i-1)}
$$
instead. By the definitions of $G$ and $t'_i$ preceding Lemma \ref {th4.4} and the remark there, it is not difficult to show that $\#{\mathcal B} = \#{\mathcal B}'$ again, hence the assertion of the lemma also hold in this case.   \qquad $\Box$

\bigskip

\noindent {\bf Proof of Theorem \ref {th3.6}.}  Let
$
{\mathcal D}'={\mathcal E}_0\oplus b^{\ell_1}{\mathcal E}_1\oplus\dotsm\oplus b^{\ell_k}{\mathcal E}_k.
$
Then ${\mathcal D}'$ is a strict product-form. By Theorem \ref{th3.3}, $K(b, {\mathcal D}')$ is a spectral set.
Let $\Lambda$ be a spectrum of it. Then, by Proposition \ref{th2.4}, $\Lambda$ has a form $\Lambda= \Gamma \oplus {\mathbb Z}$ with $\#\Gamma = |K(b,{\mathcal D}')|$. Define $K_{{\mathcal D}'}(x)=K_{{\mathcal D}'}^{(k)}(x)$ as in (\ref{eq3.11}). It follows that $P_{{\mathcal D}'}(x)=K_{{\mathcal D}'}(x)$.

\medskip
For the modulo product-form ${\mathcal D}$, we have $K_{{\mathcal D}}(x)\mid P_{{\mathcal D}}(x)$, and  $K_{{\mathcal D}}(x)= K_{\D'}(x)$ following from the definition. Therefore $P_{{\mathcal D}}(x)=P_{{\mathcal D}'}(x)Q(x)$ for some integral polynomial $Q(x)$. By (\ref {eq3.2}), we have
\begin{eqnarray*}
\widehat\mu_{b,{\mathcal D}}(x)
=\prod_{j=1}^{\infty}(b^{-1}
P_{{\mathcal D}}(e^{-2\pi ib^{-j}x})), \quad \widehat\mu_{b,{\mathcal D}^{\prime}}(x)
=\prod_{j=1}^{\infty}(b^{-1}
P_{{\mathcal D}^{\prime}}(e^{-2\pi ib^{-j}x})).
\end{eqnarray*}
Hence, $\widehat\mu_{b,{\mathcal D}'}$ is a factor of $\widehat\mu_{b,{\mathcal D}}$. This implies that $\Lambda$ is an orthonormal family for $L^2(\mu_{b,{\mathcal D}})$.  It follows from Lemma \ref{th4.5} that $|K(b,{\mathcal D})|=| K(b,{\mathcal D}^{\prime})|=\#\Gamma$. By Proposition \ref{th2.3}, $\Lambda$ is also a spectrum for $\mu_{b,{\mathcal D}}={\mathcal L}|_{K(b,{\mathcal D})}$.
\qquad $\Box$

\bigskip

\section{\bf Spectral property with ${\mathcal D} =\{0,a, b, c \}$}

\bigskip

 It is proved in Example \ref{th3.8} that if ${\mathcal D}$ is a tile digit set with respect to $4$, then $K(4, {\mathcal D})$ is a spectral set.  In this section, we will show that the converse also holds.
We assume that  ${\mathcal D}=\{0, a, b, c\}\subset\Z^+$ with $\text{gcd}(a,b,c)=1$. Let
$P_{{\mathcal D}}(x)=1+x^a+x^b+x^c$ be the mask polynomial of ${\mathcal D}$. Let $\mu$ be the self-similar measure associate with $(4, {\mathcal D})$, then by (\ref{eq3.2}),
\begin{equation} \label {eq5.1}
\widehat \mu (\xi) = 4^{-1}P_{\mathcal D} (e^{-2\pi i\xi/4})\widehat \mu (\xi/4) = \prod_{n=1}^\infty  (4^{-1}P_{\mathcal D} ( e^{-2\pi i\xi/4^n})).
\end{equation}
Our main theorem is

\medskip

\begin {theorem} \label {th5.1}
Let  ${\mathcal D} \subset{\mathbb Z}^+$ be a digit set with $\#{\mathcal D} = 4$ and gcd $({\mathcal D})=1$. If $\mu_{4,{\mathcal D}}$ is a spectral measure, then  ${\mathcal D} = \{0, a, 2^{t} \ell,  a +2^{t} \ell'\}$, where  $a, t, \ell, \ell'$ are odd integers.
\end{theorem}

\bigskip
We will prove Theorem \ref {th5.1} through a series of lemmas. As a direct consequence, we have

\begin {theorem} \label {th5.2}
Let  ${\mathcal D} \subset{\mathbb Z}^+$ be a digit set with $\#{\mathcal D} = 4$ and gcd $({\mathcal D})=1$. Then $K(4, {\mathcal D})$ is a self-affine tile  if and only if it is a spectral set.  In this case,  $K(4, {\mathcal D})$ has a spectrum $\Lambda={\Bbb Z} +\sum_{j=1}^k\frac{1}{4^j}\{0,1\}$ with $k = (t-1)/2$, where  $t$ is an odd integer as in Theorem \ref{th5.1} ($\Lambda={\Bbb Z}$ if $t=1$).
\end{theorem}

\medskip

\noindent {\bf Proof}  The first part of the theorem follows from Example \ref {th3.8} and Theorem \ref{th5.1}. In this case, the self-affine tile $K = K(4, {\mathcal D})$ and the self-affine measure $\mu$ are both determined by the same class of digit sets ${\mathcal D}$, and $\frac {d\mu}{dx} = \chi_K$.

 By Example \ref{th3.8}, we see that ${\mathcal D}$ is a modulo product-form defined by $\E = \{0, 1\}\oplus \{0, 2\}$. By (\ref{eq3.5}), the corresponding strict product-form $\D'$ gives a tile
$$
K(4, {\mathcal D}^{\prime})=[0, 1]\oplus \bigoplus_{j=0}^{k-1}4^j\{0, 2\} \quad  \hbox {with} \quad k=(t-1)/2.
 $$
If $t=1$, then $K(4, \D' ) = [0,1]$ so that $\Lambda = {\Bbb Z}$.  It $t>1$, then note that $\frac{1}{4}\{0, 1\}$ is a spectrum for the set $\{0, 2\}$,  hence $\frac{1}{4}\{0, 1\}\oplus \frac{1}{4^2}\{0, 1\}$ is a spectrum for $\{0, 2\}\oplus 4\{0, 2\}$. Inductively, we can show that $\sum_{j=1}^k\frac{1}{4^j}\{0, 1\}$ is a spectrum for the set $\bigoplus_{j=0}^{k-1}4^j\{0, 2\}$.
Hence, ${\Bbb Z}\oplus \sum_{j=1}^k\frac{1}{4^j}\{0, 1\}$ is a spectrum for $K(4, {\mathcal D}^{\prime})$ by Theorem \ref{th3.4},  and it is also a spectrum for $K(4, {\mathcal D})$ by Theorem \ref{th3.6}.
\qquad $\Box$

\bigskip

To prove Theorem \ref{th5.1}, we need to make a detailed analysis on the zeros of $P_{\D}(e^{2\pi i\xi})$ and $\widehat \mu(\xi)$. Let ${\Bbb O}$ denote the set of odd integers and ${\Bbb E}$ the set of even integers (including $0$).

\bigskip

\begin{Lem}\label{th5.3} The equation $P_{\mathcal D}(e^{2\pi i \xi})=0$ has a solution if and only if
two of $a, b, c $ are odd and one is even.
\end{Lem}

\medskip

\noindent {\bf Proof.}
That $P_{\mathcal D}(e^{2\pi i \xi})=0$ yields $ |1+e^{2\pi
ia\xi}|=|1+e^{2\pi i(c-b)\xi}|$. By symmetry, we see that the zeros satisfy one of the following sets of equations:
\begin{eqnarray} \label {eq5.2}
\begin{cases}
2a\xi=k\\
2(c-b)\xi=\ell,
\end{cases}
\begin{cases}
2b\xi=k\\
2(c-a)\xi=\ell,
\end{cases}
\begin{cases}
2c\xi=k\\
2(b-a)\xi=\ell,
\end{cases}
k,\ell\in {\Bbb O}.
\end{eqnarray}
The first set of equation yields
$\xi=\frac{k}{2a}=\frac{\ell}{2(c-b)}$.
It implies
$
k =\frac {as}{\text{gcd}(a,c-b)}$ and $\ell=\frac {(c-b)s^{\prime}}{\text{gcd}(a,c-b)}$,
where $s, s^{\prime}\in {\Bbb O}$.
Hence, we have
$
\frac {a(c-b)}{(\text{gcd}(a,c-b))^2}\in {\Bbb O}.
$
 This together with gcd$(a, b, c)  =1$ implies that one of the $a, b, c$ must be even and the other two are odds, and the solution $\xi$ is given by
$
\xi = \frac{j}{2 \text {gcd}(a,c-b)},\ j \in {\Bbb O} .
$
The same argument applies to the other two sets of equations.
\qquad $\Box$

\bigskip

In view of Lemma \ref{th5.3},  we assume in the rest of this section, that $b\in {\Bbb E}$, $a<c$  and $a,c\in {\Bbb O}$.
Let $b=2^t\ell$ and $c-a=2^{t^{\prime}}\ell^{\prime}$, where $t,
t^{\prime}\ge 1$ and $\ell,\ell^{\prime}\in {\Bbb O}$.  Hence ${\mathcal D} = \{0, a, 2^t\ell, a+2^{t'}\ell'\}$. For simplicity, we let
$
p_1=\text{gcd}(a, c-b), \quad p_2=\text{gcd}(c, b-a) \ \  \hbox {and} \ \   p_3=\text{gcd}(\ell, \ell').
$
Note that $p_i  \in {\Bbb O}$ and gcd$(p_i, p_j) =1$ for $i \not =j$. Let
$$
{\mathcal R}_i=\frac 1{2p_i}{\Bbb O} , \ \ i=1,2 \quad \hbox {and} \quad {\mathcal R}_3=\frac 1{2^{t+1}p_3}{\Bbb O}.
$$
We observe that ${\mathcal R}_1$ and ${\mathcal R}_2$ correspond to the roots of the first and the third sets of equations in (\ref{eq5.2}), and ${\mathcal R}_3$ corresponds to the second set of equation.  For ${\mathcal R}_3$, we note that the solution must satisfy ${b(c-a)}/{(\text{gcd}(b,c-a))^2}\in {\Bbb O}$ (see the proof of Lemma \ref{th5.3}). This implies that $t=t'$. Hence, if  $t\not =t'$, then there is no solution belonging to
${\mathcal R}_3$. Summarizing, we have

\medskip

\begin{Lem} \label{th5.4}
 Let ${\mathcal D} = \{0, a, 2^t\ell, a+2^{t'}\ell'\}$ where $t,
t^{\prime}\ge 1$ and $\ell,\ell^{\prime}\in {\Bbb O}$.  Let
${\mathcal R}(P_{{\mathcal D}})$ be the roots of $P_{\mathcal
D}(e^{2\pi i \xi})$. Then 
$$ {\mathcal R}(P_{{\mathcal
D}}) =\left\{
                                 \begin{array}{ll}
                                   {\mathcal R}_1\cup {\mathcal R}_2 \cup {\mathcal R}_3, & \hbox{if $t=t'$;} \\
                                   {\mathcal R}_1\cup {\mathcal R}_2, & \hbox{if $t\ne t'$.}
                                 \end{array}
                               \right.
$$
\end{Lem}

\bigskip
To consider the zeros for $\widehat{\mu}(\xi)$, we let  ${\mathcal Z}(\widehat{\mu})=\{\xi\in{\mathbb R}:\widehat{\mu}(\xi)=0\}$. Let
\begin{eqnarray}\label{eq5.3}
{\mathcal Z}_i={\bigcup}_{k=1}^\infty 4^k{\mathcal R}_i, \quad \mbox{for} \ \ i=1, 2, 3.
\end{eqnarray}
Then it follows from (\ref {eq5.1}) that
$$
{\mathcal Z}(\widehat{\mu})\ =\  {\bigcup}_{k=1}^\infty 4^k{\mathcal
R}(P_{{\mathcal D}}) =\left\{
                                      \begin{array}{ll}
                                        {\mathcal Z}_1\cup{\mathcal
Z}_2\cup{\mathcal Z}_3, & \hbox{if  $t=t'$;} \\
                                        {\mathcal Z}_1\cup{\mathcal
Z}_2, & \hbox{if $t\ne t'$.}
                                      \end{array}
                                    \right.
$$

\medskip

\begin{Lem}\label{th5.5}
Let $0\in\Lambda\subset\Bbb R$ be an orthogonal set for $\mu$. Then $\Lambda\setminus\{0\}\subseteq {\mathcal Z}_i\cup {\mathcal Z}_3$ for $i=1$ or $2$.
\end{Lem}

\medskip

\noindent {\bf Proof.} By the orthogonal property,  $\Lambda -\Lambda \subset {\mathcal Z}(\widehat{\mu})$.  Hence, it suffices to prove the lemma by showing that for $\lambda_1 \in \Lambda\cap ({\mathcal Z}_1\setminus {\mathcal Z}_2)$ and $\lambda_2 \in \Lambda\cap ({\mathcal Z}_2\setminus {\mathcal Z}_1)$, then $\lambda_1 - \lambda_2 \not \in {\mathcal Z}(\widehat{\mu})$.
Otherwise, let
\begin{eqnarray*}
\lambda_1 = 4^{k_1}\frac {a_1}{2p_1}\in \Lambda\cap({\mathcal Z}_1\setminus {\mathcal Z}_2) \quad {\hbox{and}} \quad \lambda_2=4^{k_2}\frac {a_2}{2p_2}\in
 \Lambda\cap({\mathcal Z}_2\setminus {\mathcal Z}_1).
\end{eqnarray*}
Then $p_i \nmid a_i$  for $a_i\in {\Bbb O}$ and $i=1,2$. By the orthogonal property, we have
either $\lambda_1-\lambda_2=4^k\frac{a}{2p}$ where $p=p_1$ or $p_2$, or $\lambda_1-\lambda_2= 4^k\frac{a}{2^{1+t}p_3}$ for some
$a\in {\Bbb O}$.
For the first case, we obtain that
$$
4^{k_1}\frac {a_1}{2p_1}-4^{k_2}\frac {a_2}{2p_2}=4^{k}\frac {a}{2p}.
$$
Without loss of generality, we assume that $p=p_1$. Then
by $\text{gcd}(p_1, p_2)=1$, it is easy to see that  $p_2\mid a_2$, which is a contradiction. For the second case, we have
\begin{equation} \label {eq5.4}
4^{k_1}\frac {a_1}{2p_1}-4^{k_2}\frac {a_2}{2p_2}=4^k\frac{a_3}{2^{1+t}p_3}.
\end{equation}
Since $\text{gcd}(p_i, p_j)=1$ for $ i\ne j$, we have $p_i\mid a_i$ for $i=1,2,3$. This yields a contradiction again.
\qquad $\Box$

\bigskip

We can improve the above lemma further to smaller subsets of ${\mathcal Z}_i \cup {\mathcal Z}_3, \ i = 1$ or $2$.  Let
$$
{\mathcal Z}_0 = \bigcup_{k=1}^\infty \frac {4^k}{2} {\Bbb O} \ \ (\subset {\mathcal Z}_1 \cap {\mathcal Z}_2), \qquad {\mathcal Z}_3' = \bigcup_{k=1}^{\infty}\frac {4^k}{2^{t+1}}{\Bbb O} \ \ (\subset {\mathcal Z}_3).
$$

\medskip

\begin{Lem}\label{th5.6}
Let $0\in\Lambda\subset\Bbb R$ be an orthogonal set of  $\mu$. Then
 either $\Lambda \setminus \{0\} \subseteq {\mathcal Z}_0 \cup {\mathcal Z}_3$\ \ or \ $\Lambda \setminus \{0\}\subseteq {\mathcal Z}_j \cup {\mathcal Z}'_3, \ j =1 \ \hbox {or} \ 2$. In particular, if $t\in {\Bbb E}$, then we have $\Lambda \setminus
\{0\}\subseteq {\mathcal Z}_3$ \ or \
$
\Lambda \setminus \{0\}\subseteq {\mathcal Z}_j\cup {\mathcal Z}'_3, \ j=1\, \ \mbox{or}\, \ 2.
$
\end{Lem}

\medskip

\noindent {Proof.}  According to Lemma \ref{th5.5},  we can assume, without loss of generality,  that $\Lambda \subseteq {\mathcal Z}_1\cup{\mathcal Z}_3$. We claim that it is impossible to have
$$
\Lambda\cap ({\mathcal Z}_1\setminus {\mathcal Z}_0)\ne \emptyset \quad \mbox{and} \quad
\Lambda\cap ({\mathcal Z}_3\setminus{\mathcal Z}'_3)\ne \emptyset.
$$
Suppose, otherwise, that $\lambda_1$ and $\lambda_2$ are in the two sets respectively, then
$\lambda_1=\frac{4^{k_1}a_1}{2p_1}$ with $p_1\nmid a_1$ and $\lambda_2=\frac{4^{k_2}a_2}{2^{1+t}p_3}$ with $p_3\nmid a_2$,
where $a_i\in {\Bbb O}$, $i=1, 2$. Since $\lambda_1-\lambda_2\in {\mathcal Z}_1\cup{\mathcal Z}_2\cup{\mathcal Z}_3$, by \eqref{eq5.4}, we see that
 $\lambda_1-\lambda_2\not\in {\mathcal Z}_2$. If $\lambda_1-\lambda_2\in{\mathcal Z}_1$, we have
$$
\frac{4^{k_1}a_1}{2p_1}-\frac{4^{k_2}a_2}{2^{1+t}p_3}
=\frac{4^{k}a}{2p_1}
 $$
for  some $ a\in {\Bbb O}$, which  implies $p_3|a_2$,  and is a contradiction. If $\lambda_1-\lambda_2\in{\mathcal Z}_3$, we have
\begin{eqnarray}\label{eq5.5}
\frac{4^{k_1}a_1}{2p_1}-\frac{4^{k_2}a_2}{2^{1+t}p_3}
=\frac{4^{k}a}{2^{1+t}p_3} \quad {\hbox {for \ some} } \ \ a\in {\Bbb O}.
\end{eqnarray}
This forces $p_1| a_1$, which contradicts with our assumption again. Hence the claim holds.
If $\Lambda\cap ({\mathcal Z}_1\setminus{\mathcal Z}_0)=\emptyset$, then $\Lambda \setminus \{0\} \subseteq {\mathcal Z}_0
 \cup {\mathcal Z}_3$. On the other hand, if $\Lambda\cap\left({\mathcal Z}_1\setminus{\mathcal Z}_0\right)\ne \emptyset$, then by the claim, $\Lambda\cap\left({\mathcal Z}_3\setminus{\mathcal Z}_3'\right)=\emptyset$.
Therefore we have
$\Lambda \setminus \{0\}\subseteq {\mathcal Z}_1 \cup {\mathcal Z}_3'.$

 Particularly, if $t\in\Bbb E$, then ${\mathcal Z}_3
 \cup {\mathcal Z}_0 = {\mathcal Z}_3$ since $\frac {4^k}2a=
 \frac {4^{k+t/2}}{2^{1+t}}a\in{\mathcal Z}_3$ for any $a\in {\Bbb O}$.
\qquad $\Box$

\bigskip

In the following we will  exclude the cases $t\in {\Bbb E}$ (in Lemma \ref{th5.6}) and $t\not = t'$ to conclude Theorem \ref {th5.1}.

\medskip

Suppose  $b=2^t\ell$ for some $t\in {\Bbb E}$ and $\ell\in {\Bbb O}$. In view of Lemma \ref {th5.6}, we first consider  $\Lambda\setminus\{0\}\subseteq {\mathcal Z}_3$. Let us write
$\Lambda=\{\lambda_n\}_{n=0}^\infty := \{\frac {a_n} {2^{1+t}p_3}\}_{n=0}^\infty$, where $a_0=0$ and $a_n\in \{4^kd: k\ge 1, d\in{\Bbb O}\}$ for $n\ge1$. Hence,
for each $n\ge 1$, $a_n$ can be expressed uniquely by
\begin{eqnarray}\label{eq5.6}
a_n=l_{n, 1}4+l_{n, 2}4^2+\cdots+l_{n, \alpha(n)}4^{\alpha(n)}=\sum_{k=1}^\infty l_{n, k}4^k,
\end{eqnarray}
where all $l_{n, k}\in\{-1, 0, 1, 2\}$, $l_{n, \alpha(n)}\ne 0$, $l_{n, k}=0$ for $k>\alpha(n)$.
 (We use $\{-1, 0, 1, 2\}$ in the expansion instead of the $\{0,1, 2,3\}$ so as to cover all the $a_n$ positive or negative.)
The orthogonal property implies that, for any $n\ne n'\in\Z^+$, $
\frac 1{2^{1+t}p_3} (a_n -a_{n'}) \in {\mathcal Z}(\widehat{\mu})
\subseteq {\mathcal Z}_1\cup{\mathcal
Z}_2\cup{\mathcal Z}_3$. By $t\in\Bbb E$ and the same argument as in
\eqref{eq5.5}, we have
\begin{eqnarray}\label{eq5.7}
a_n-a_{n'}\in \{4^kd: k\ge 1, d\in{\Bbb O}\}.
\end{eqnarray}
Let $s$ be the first index such that $l_{n, s}\ne l_{n', s}$. It follows from \eqref{eq5.7} that one of $l_{n, s}$ and $l_{n', s}$ is odd and the other one
 is even. Thus either $l_{n, s}\in\{-1, 1\}$ or $l_{n, s}\in\{0, 2\}$ (we remark that the orthogonality implies that $l_{n, k}$ depends on all $l_{n, i}$ ($1\le i<k$)). Hence, for $k\ge 1$, the set
$$
\A_k:=\{l_{n,1}4+l_{n,2}4^2+\cdots+l_{n,k}4^k: n\ge 1\}
$$
 has at most $2^k$ elements.

\bigskip

\begin{Lem}\label{th5.7}
Let $b=2^t\ell$ for some $t\in\Bbb E$ and $\ell\in\Bbb O$, then any  orthogonal set $\Lambda \subset {\mathcal Z}_3\cup \{0\}$  cannot be a spectrum of $\mu$.
\end{Lem}

\medskip

\noindent {\bf Proof.}   We will show that  $Q(\xi)=\sum_{n=0}^\infty |\widehat{\mu}(\xi+ {\lambda_n})|^2 \not \equiv 1$. Let $m = t/2$ and let
$\nu_m=\delta_{4^{-1}{\mathcal D}}\ast\delta_{4^{-2}{\mathcal D}}\ast\cdots\ast\delta_{4^{-m}{\mathcal D}}$. Then $\widehat \nu_m = \prod_{n=1}^m 4^{-1} P_{\mathcal D}(e^{-2\pi i \xi/4^n})$. By \eqref{eq5.1},
we can write
\begin{eqnarray*}Q(\xi)=\sum_{n=0}^\infty \Big|\widehat{\nu}_{m}\big (\xi+ \lambda_n \big )\ \widehat{\mu}\big (4^{-m}
(\xi+ \lambda_n)\big )\Big |^2.
\end{eqnarray*}
Write
 $a_n=a_{n}'+a_{n}^{\prime \prime}$, where
$a_{n}'=\sum_{i=1}^{m}l_{n,i}4^i \ (\in \A_m)$, and correspondingly,  $\lambda_n = \lambda_n' + \lambda_n^{\prime \prime}$.  For $\gamma \in \A_m$, let
$
\Lambda_\gamma = \{\lambda_n:  a'_n= \gamma\}$ and $\lambda_\gamma = \frac \gamma{2^{1+t}p_3}.$
 By the  $4^m$-period of $\widehat{\nu}_{m}(\xi)$, we can find
$0\leq r(\xi, \gamma)< 4^m$ such that
$$
\big |\widehat{\nu}_{m}(\xi+\lambda_\gamma +r(\xi, \gamma))\big |\ =\ \max\big \{\big |\widehat{\nu}_{m}(\xi+{\lambda_n}) \big |: \lambda_n\in \Lambda_\gamma \big \}.
$$
Also it is clear that $4^{-m}\Lambda_\gamma$ is an orthogonal set for the measure  $\mu(4^m\ \cdot)$. Hence by (\ref{eq2.2}),
$\sum_{\{n: \lambda'_n=\lambda_\gamma\}}\big | \widehat{\mu}(4^{-m}(\xi+ {\lambda_n}))\big |^2\le 1$. It follows that
\begin{eqnarray}\label {eq5.8}
Q(\xi)
 & \le & \sum_{{\gamma} \in {\mathcal A}_m}\Big |\widehat{\nu}_{m}(\xi+\lambda_\gamma + r(\xi, \gamma)/2^{1+t}p_3) \Big |^2 \Big ( \sum_{\{n: \lambda_n' = \lambda_\gamma\}}\Big | \widehat{\mu}(4^{-m}(\xi+ \lambda_n))\Big |^2\Big )\nonumber\\
 &\leq & \sum_{{\gamma} \in {\mathcal A}_m}\Big |\widehat{\nu}_{m}(\xi+\lambda_ \gamma + r(\xi, \gamma)) \Big |^2.
\end{eqnarray}
Note that $\Lambda_m =\{ \lambda_\gamma: \gamma \in \A_m \}$ is an orthogonal set of $\nu_m$. However, it cannot be a spectrum of $\nu_m$ for   $m$ large enough, in view of the fact that $\#\A_m \le 2^m$,  and  $L^2(\nu_m)$ has dimension $4^m$. Hence it follows from (\ref {eq2.2}) and (\ref{eq5.8}) that $Q(\xi)<1$.  This implies that $\Lambda$ cannot be a spectral set, and completes the proof. \qquad $\Box$

\bigskip

\begin{Lem}\label{th5.8}
Let $b=2^t\ell$ for some $t\in\Bbb E$ and $\ell\in\Bbb O$. Then any  orthogonal set $\Lambda \subset {\mathcal Z}_i \cup {\mathcal Z}_3'\cup\{0\}, \ i =1$ or $2$,   cannot be a spectrum of $\mu$.
\end{Lem}

\medskip

\noindent {\bf Proof.}  Let us consider $i=1$.
Note that for $\frac{4^ka}{2^{t+1}}\in {\mathcal Z}_3'$ with $k>t/2$ and $a\in\Bbb O$, we have $\frac{4^ka}{2^{1+t}}=\frac{4^{k-\frac{t}{2}}a}{2}\in{\mathcal Z}_1$; also for $\frac{4^ka}{2^{t+1}}\in {\mathcal Z}_3'$ with $k<t/2$, we have $\frac{4^ka}{2^{t+1}}\in \frac{1}{2^t}{\mathcal Z}_1$. These imply that
${\mathcal Z}_1\cup {\mathcal Z}_3'\subset \frac{1}{2^t}{\mathcal Z}_1$.
The lemma now follows by observing that the proof of  Lemma \ref{th5.7} also applies to the statement: any  orthogonal set
$\Lambda \setminus\{0\}\subset \frac{1}{2^t}{\mathcal Z}_1=\frac {1}{2^{t+1}p_1}\bigcup_{k=1}^{\infty}4^k{\Bbb O}$ cannot be a spectrum of $\mu$.  \qquad $\Box$

\bigskip

\begin {Lem} \label {th5.9}
For $t\ne t'$, then any orthogonal set $\Lambda$ of $\mu$  cannot be a spectrum of $\mu$.
\end{Lem}

\medskip

\noindent {\bf Proof.} By Lemma \ref {th5.4} and Lemma \ref{th5.5}, we see that for $t \ne t'$, then $\Lambda \setminus \{0\} \subset {\mathcal Z}_i$ for $i =1,2$. As ${\mathcal Z}_i \subset \frac 1{2^s}{\mathcal Z}_i$ for any given  $s\in\Bbb E$,  the lemma follows by the last statement in the proof of Lemma \ref {th5.8} \qquad $\Box$

\bigskip

\noindent {\bf Proof of Theorem \ref{th5.1}} In Lemma \ref {th5.7} and Lemma \ref{th5.8}, we have excluded the case $t\in {\Bbb E}$, and in Lemma \ref {th5.9} the case $t\not = t'$. Hence the remaining case is ${\mathcal D} = \{0, a, 2^t\ell, a+2^{t'}\ell'\}$ with $t=t^{\prime}\in\Bbb O$ as in the theorem.  \qquad $\Box$

\bigskip

\bigskip

\section{\bf Remarks}

 For the main theorems in Sections 4 and 5, the modulo product-forms  and the product-forms  under consideration rely on the ``strict" product $\E =
\E_0 \oplus \cdots \oplus\E_k = {\Bbb Z}_b$.
However, not all modulo product-forms  can be constructed through such strict product (the cases where $b=p^{\alpha}$ or $pq$ in Proposition \ref {th3.7} are exceptional).
 In fact, the example  appeared in [PW] (see also [dB]) with $b=72 (= 2^33^2)$ and
$$
{\mathcal A}=\{0, 8, 16, 18, 26, 34\} \ \ \ \ \ \ \ {\mathcal B}=\{0, 5, 6, 9, 12, 29, 33, 36, 42, 48, 53, 57\}
$$
satisfying that
${\mathcal A}\oplus {\mathcal B}\equiv {\mathbb Z}_{72} (\text{mod} \ 72)$, we can show that $\D = {\mathcal A}\oplus {\mathcal B}$ cannot be expressed as a  modulo product-form defined by a strict product-form. If we assume that ${\mathcal A}\oplus {\mathcal B}$ is a modulo product-form defined by a strict product-form, according to Definition \ref{th3.5}, its corresponding strict product-form is ${\mathcal D}^{\prime}={\mathbb Z}_{72}$. In all irreducible decompositions of ${\mathbb Z}_{72}$ (see (\ref{eq4.1})), either $2^3\cdot 3 \ {\mathbb Z}_3$ or $2^2\cdot 3^2 \ {\mathbb Z}_2$ must be an irreducible component of ${\mathbb Z}_{72}$. Then, either $\Phi_{2^{3}}(x^{3^2})$ or $\Phi_{3^{2}}(x^{2^3})$ must be a factor of $P_{{\mathcal D}^{\prime}}(x)=K_{{\mathcal D}^{\prime}}(x)$, which is defined in Definition \ref{th3.5}.  Since $K_{{\mathcal D}}(x)=K_{{\mathcal D}^{\prime}}(x)$ and $K_{{\mathcal D}}(x)\mid P_{{\mathcal D}}(x)=P_{{\mathcal A}}(x) P_{{\mathcal B}}(x)$, then either $\Phi_{2^{3}}(x^{3^2})$ or $\Phi_{3^{2}}(x^{2^3})$ must be a factor of $P_{{\mathcal A}}(x)$ or $P_{{\mathcal B}}(x)$.
This is not impossible for our ${\mathcal A}$ and ${\mathcal B}$ since
$$
P_{{\mathcal A}}(x)=\Phi_3(x^{2^3})\Phi_{2^2}(x^{3^2}), \ P_{{\mathcal B}}(x)=\Phi_{2^3}(x)\Phi_{3^2}(x)\Phi_{18}(x)\Phi_{72}(x)(\Phi_{24}(x)+x^5\Phi_{36}(x)).
$$
 This yields a contradiction.  A challenging question is

\bigskip

 {\bf Q1}. Show that the product-form and the modulo product-form tile defined by $\E =
\E_0 \oplus \cdots \E_k \equiv {\Bbb Z}_b \ (\hbox {mod} \ b)$ are  spectral sets.

\medskip

We note that in Proposition \ref{th3.2}, the product-form tile satisfies $K(b,\D) = K(b, \E) +{\mathcal A}$ with ${\mathcal A}$ an integer tile. It will also be interesting to study the Fuglede conjecture for this finite set ${\mathcal A}$, which will yield the spectral property for $K(b,\D)$. There are investigations for such case. It is known that an integer tile ${\mathcal A}$ is equivalent to spectral set if $\#{\mathcal A}$ is the product of two prime powers ([CM], [$\L$2]), but the general case is still unsolved.

\bigskip

The converse of the above problem is more intricate, we formulate the question as

\medskip

{\bf Q2}. Let $\mu_{b, \D}$ be a self-similar measure (with $b=\#\D$) and is spectral, does it imply $K(b, \D)$ is a tile?

\medskip

We have seen in Section 5 that the statement is true, but the proof is indirect, and is difficult to extend. It will be instructive to obtain another more direct proof.  In particular, it is worthwhile to study  the cases  for  $b=p^{\alpha}$ or $b=p^\alpha q^\beta$, where $p, \ q$ are distinct primes and $\alpha, \beta \ge 1$,  as there are more results on the geometry of numbers we can use ([dB],[CM], [LLR1,2]).

\medskip
So far we have not considered the spectral problem other than dimension one. For the product-forms, it can be defined on higher dimension. The modulo product-form poses more question, as it is defined through the cyclotomic polynomial, which is not so clear how to extend to the higher dimension.


\bigskip
\bigskip

\begin{thebibliography}{9999}

\bibitem [B] {B}
{\sc C. Bandt}, {\it Self-similar sets V. Integer matrices and fractal tilings of ${\mathbb R}^n$},
 Proc. Amer. Math. Soc. \textbf{112}(1991), 549--562



\bibitem [dB] {dB}
{\sc D. de Bruijn}, {\it On the factorization of cyclic group},
 Indag. Math. \textbf{17}(1955), 370--377


%\bibitem [BM] {BM}
%{\sc D. Bose, S. Madan}, {\it Spectrum is periodic for $n$-intervals.}
%\newblock J. Funct. Anal. \textbf{260}, 308--325
 % (2011)

\bibitem [CM] {CM}
{\sc E. Coven and A. Meyerowitz}, {\it Tiling the integers with translates of one finite set},
 J. Algebra \textbf{212} (1999), 161--174


\bibitem [D] {D}
{\sc X.-R. Dai,} {\it When does a Bernoulli convolution admit a
spectrum?}
 Adv. Math. \textbf{231}(2012), 187--208


\bibitem [DHL] {DHL}
{\sc X.-R. Dai, X.-G. He and C.-K. Lai,} {\it Spectral structure of
Cantor measures with consecutive digits}, Adv. Math.
\textbf{242}(2013), 1681--1693



\bibitem [DHS] {DHS}
{\sc D. Dutkay, D.-G. Han and Q.-Y. Sun,} {\it On the spectra of a
Cantor measure}, Adv. Math. \textbf{221}(2009), 251--276


%\bibitem [DL1] {DL1}
%{\sc D. E. Dutkay, C. K. Lai,} {\it Uniformity of measures with Fourier frames.}
%\newblock arXiv:1202.6028v1 \textbf{online}, 1--25
 % (2012)

%\bibitem [DL] {DL}
%{\sc D. Dutkay and C. K. Lai,} {\it Some reductions of the spectral set conjecture to integers}, preprint



\bibitem [DJ1]{DJ1}
{\sc  D.  Dutkay and P.  Jorgensen,} {\it Quasiperiodic spectra and orthogonality for iterated function system measures},
 Math. Z. \textbf{261}(2009), 373--397


%\bibitem [DJ2]{DJ2}
%{\sc D.  Dutkay, P. Jorgensen,} {\it Fourier duality for fractal measures with affine scales.}
%\newblock Math. Comp. \textbf{81}, 2253--2273
  %(2012)


\bibitem [DJ2] {DJ2}
{\sc D. Dutkay and P.  Jorgensen,} {\it On the universal tiling conjecture in dimension one.}
\newblock J. Fourier Anal. Appl. \textbf{19}(2013), 467--477



\bibitem [F] {F}
{\sc B. Fuglede,} {\it Commuting self-adjoint partial differential operators and a group theoretic problem},
\ J. Funct. Anal. \textbf{16}(1974), 101--121



\bibitem [H] {H}
{\sc J. Hutchinson,} {\it Fractals and self-similarity},
Indiana Univ. Math. J. \textbf{30}(1981), 713--747

\bibitem [HL]{HL}
{\sc X.-G. He, K.-S. Lau,} {\it Characterization of tile digit sets
with prime determinants.}
\newblock Appl. Comput. Harmon. Anal.  \textbf{16}(2004), 159--173


\bibitem [HLL] {HLL}
{\sc X.-G. He, C.-K. Lai and K.-S. Lau,} {\it Exponential spectra in
$L^2(\mu)$}, Appl. Comput. Harmon. Anal.  {\bf 34}(2013), 327-338

\bibitem [HuL]{HuL}
{\sc T.-Y. Hu and K.-S. Lau,} {\it Spectral property of the
Bernoulli convolutions},
 Adv. Math. \textbf{219}(2008), 554--567

\bibitem [IKP] {IKP}
{\sc  A. Iosevich, N. Katz and S. Pedersen,} {\it Fourier bases and a distance problem of Erd\"os},
 Math. Res. Lett. \textbf{6}(1999), 251--255



\bibitem [IKT1]{IKT1}
{\sc A. Iosevich, N.  Katz and T. Tao,} {\it Convex bodies with a point of curvature do not have Fourier bases},
 Amer. J. Math. \textbf{123}(2001), 115--120


\bibitem [IKT2]{IKT2}
{\sc A. Iosevich,  N.  Katz and T. Tao,} {\it The Fuglede spectral conjecture holds for convex planar domains},
 Math. Res. Lett. \textbf{10}(2003), 559--569


%\bibitem [IK] {IK}
%{\sc A. Iosevich, M. N. Kolountzakis,} {\it Periodicity of the spectrum in dimension one.}
%\newblock arXiv:1108.5689v2 \textbf{10}, 559--569
 % (2012)

%\bibitem [JKS] {JKS}
%{\sc P. T. E. Jorgensen, K. Kornelson, K. Shuman,} {\it Families of spectral sets for Bernoulli convolutions.}
%\newblock J. Fourier Anal. Appl. \textbf{17}, 431--456
 % (2011)

\bibitem [JP] {JP}
{\sc P.  Jorgensen and S. Pedersen,} {\it Dense analytic subspaces in fractal $L^2$-spaces},
 J. Anal. Math. \textbf{75}(1998), 185--228



%\bibitem                                                           [K]{K}
%{\sc M. N. Kolountzakis,} {\it Non-symmetric convex domains have no basis of exponentials},
%\newblock      Illinois      J.      Math.      \textbf{44},     542--550
%  (2000)

\bibitem[K]{K}
{\sc R. Kenyon,} {\it Self-replicating tilings},  P. Walters (Ed.), Symbolic Dynamics and Its Applications, Contemp. Math.
\textbf{135}(1992), 239--264

\bibitem [KM1]{KM1}
{\sc M. N. Kolountzakis and M. Matolcsi,} {\it Complex Hadamard matrices and the spectral set conjecture},
 Collec. Math. \textbf{Vol. Extra}(2006), 281--291


\bibitem [KM2]{KM2}
{\sc M. N. Kolountzakis and M. Matolcsi,} {\it Tiles with no spectra},
 Forum Math. \textbf{18}(2006), 519--528


\bibitem [$\L$1]{L1}
{\sc I. $\L$aba,} {\it Fuglede's conjecture for a union of two intervals},
 Proc. Amer. Math. Soc. \textbf{129}(2001), 2965--2972


\bibitem [$\L$2]{L2}
{\sc I. $\L$aba,} {\it The spectral set conjecture and multiplicative properties of roots of polynomials},
\newblock J. London Math. Soc. \textbf{65}(2002), 661--671


\bibitem[$\L$W]{LW}
{\sc I. $\L$aba and Y. Wang,} {\it On spectral Cantor measures},
J. Funct. Anal. \textbf{193}(2002), 409--420


\bibitem[LW1]{LW1}
{\sc J. Lagarias and Y. Wang,} {\it Integral self-affine tiles in
\mbox{${\mathbb R}^n$}
  \mbox{I}. \mbox{Standard} and nonstandard digit sets},
 J. London Math. Soc. \textbf{54}(1996), 161--179


\bibitem [LW2]{LW2}
{\sc J. Lagarias and Y. Wang,} {\it Tiling the line with translates
of one tile},
 Invent. Math. \textbf{124}(1996), 341--365


\bibitem [LW3]{LW3}
{\sc J.  Lagarias and Y. Wang,} {\it Self-affine tiles in
\mbox{${\mathbb R}^n$}},
 Adv.  Math. \textbf{121}(1996), 21--49


\bibitem [LW4]{LW4}
{\sc J. C. Lagarias and Y. Wang,} {\it Spectral sets and
factorization of finite abelian groups},
 J. Funct. Anal. \textbf{145}(1997), 73--98


\bibitem [LLR1]{LLR1}
{\sc C.-K. Lai, K.-S. Lau and H. Rao,} {\it Spectral structure of
digit sets of self-similar tiles on ${\mathbb R}^1$},
 Trans. Amer. Math. Soc. \textbf{365}(2013), 3831--3850


\bibitem [LLR2]{LLR2}
{\sc C.-K. Lai, K.-S. Lau and H. Rao,} {\it Classification of tile
digit sets as product-forms}, preprint



\bibitem [LR]{LR}
 {\sc K.-S. Lau and H. Rao,} {\it On one-dimensional self-similar tilings and the $pq$-tilings},
Trans. Amer. Math. Soc. \textbf{355}(2003), 1401--1414

%
%\bibitem [M]{M}
%{\sc M. Matolcsi,} {\it Fuglede's conjecture fails in dimension 4},
%\newblock Proc. Amer. Math. Soc. \textbf{133}(2005), 3021--3026
\bibitem [O]{O}
{\sc A. M. Odlyzko,} {\it Non-negative digit sets in positional
number systems},
 Proc. London Math. Soc. \textbf{37}(1978), 213--229

\bibitem [P]{P}
{\sc S. Pedersen,} {\it Spectral sets whose spectrum is a lattice with a base},
J. Funct. Anal. \textbf{141}(1996), 496--509


\bibitem [PW]{PW}
{\sc S. Pedersen and Y. Wang,} {\it Universal spectra, universal tiling sets and the spectral set conjecture},
 Math. Scand. \textbf{88} (2001), 246--256

\bibitem [T]{T}
{\sc T. Tao,} {\it Fuglede's conjecture is false in 5 and higher
dimensions},
 Math. Res. Lett. \textbf{11}(2004), 251--258
\end{thebibliography}
\end{document}